\documentclass{amsart}
\usepackage{amssymb}
\usepackage{enumerate}
\usepackage{graphicx}
\usepackage{comment}

\def\he{\hat{e}}
\def\op{{\rm op}}

\def\np{\bigskip\noindent}
\def\nl{\smallskip\noindent}
\def\suppress#1{}
\def\cl#1{{#1}^{\rm cl}}

\def\lijntje{\vrule height2.4pt depth-2pt width0.5in}

\def\vlijntje{\vrule height0.45in depth0.4pt width0.4pt}
\def\vlijn{\buildrel {\hbox to 0pt{\hss$\textstyle\circ$\hss}}\over\vlijntje}

\def\vtriple#1\over#2\over#3{\mathrel{\mathop{\kern0pt #2}\limits_{\hbox
to 0pt{\hss$#1$\hss}}^{\hbox to 0pt{\hss$#3$\hss}}}}
\def\rvtriple#1\over#2\over#3{\mathrel{\mathop{\kern0pt #2}\limits_{\hbox
to 0pt{\hss$#3$\hss}}^{\hbox to 0pt{\hss$#1$\hss}}}}

\def\Dnrev{\vtriple{\scriptstyle n}\over\circ\over{}\kern-1pt\lijntje\kern-1pt
\vtriple{\scriptstyle{n-1}}\over\circ\over{}
\cdots\cdots\vtriple{\scriptstyle 4}\over\circ\over{}\kern-1pt\lijntje\kern-1pt
\vtriple{\scriptstyle 3}\over\circ\over{\buildrel
{\scriptstyle 2}\over\vlijn}\kern-1pt\lijntje\kern-1pt
\vtriple{\scriptstyle 1}\over\circ\over{}\kern-1pt}

\def\DnS{\vtriple{\scriptstyle 0}\over\circ\over{}\kern-1pt\phantom{\lijntje}\kern-1pt
\vtriple{\scriptstyle n-2t}\over\circ\over{}\kern-1pt\lijntje\kern-1pt
\vtriple{\scriptstyle n-2t-1}\over\circ\over{}
\cdots\cdots\vtriple{\scriptstyle 4}\over\circ\over{}\kern-1pt\lijntje\kern-1pt
\vtriple{\scriptstyle 3}\over\circ\over{\buildrel
{\scriptstyle 2}\over\vlijn}\kern-1pt\lijntje\kern-1pt
\vtriple{\scriptstyle 1}\over\circ\over{}\kern-1pt}

\def\AO{{\mathcal  A}}
\def\MY{{\mathcal  Y}}

\def\Supp{{\rm Supp}}
\def\M{M}
\def\A{{\rm A}}
\def\rr{{b}}

\def\TL{{\rm {\bf TL}}}

\def\KT{{\rm {\bf KT}}}
\def\BMW{{\rm {\rm B}}}
\def\Br{{\rm {\rm Br}}}

\def\BrM{{\rm BrM}}

\def\ADE{{\rm ADE}}
\def\D{{\rm D}}
\def\E{{\rm E}}

\def\np{\medskip}
\def\nl{\smallskip\noindent}

\def\Proj{{\rm Proj}}

\def\und{\underline}

\newtheorem{Thm}{\bf{Theorem}}[section]
\newtheorem{Def}[Thm]{\bf{Definition}}
\newtheorem{Defs}[Thm]{\bf{Definitions}}

\newtheorem{Lm}[Thm]{\bf{Lemma}}

\newtheorem{Prop}[Thm]{\bf{Proposition}}
\newtheorem{Cor}[Thm]{\bf{Corollary}}
\newtheorem{Remark}[Thm]{\bf{Remark}}

\newtheorem{Notation}[Thm]{Notation}
\newtheorem{Alg} [Thm]{\bf{Algorithm}}

\renewcommand{\phi}{\varphi}

\DeclareMathOperator{\het}{ht}

\def\Z{{\mathbb Z}}
\def\Q{{\mathbb Q}}
\def\R{{\mathbb R}}

\def\alg{{A}}

\newcommand{\ov}{\overline}

\newcommand{\N}{{\mathbb{N}}}

\newcommand{\eps}{\varepsilon}

\def\a{\alpha}
\def\b{\beta}
\def\c{\gamma}

\def\homog{\rightsquigarrow}
\def\isog{\leftrightsquigarrow}

\def\CoxDiag{M}

\author{Arjeh M. Cohen
\& Di\'{e} A.H. Gijsbers
\& David B. Wales}
\address{Arjeh M. Cohen\\
Department of Mathematics and Computer Science\\
Eindhoven University of Technology\\
POBox 513\\
5600 MB Eindhoven\\
The Netherlands}
\email{A.M.Cohen@tue.nl}
\address{Di\'{e} A.H. Gijsbers}
\email{dahgijsbers@gmail.com}
\address{David B. Wales\\
Mathematics Department\\
Sloan Lab\\
Caltech\\
Pasadena, CA 91125\\
USA}
\email{dbw@its.caltech.edu}

\title{The Birman--Murakami--Wenzl Algebras of Type $\D_n$}

\date{\today}

\begin{document}

\begin{abstract}
The Birman--Murakami--Wenzl algebra (BMW algebra) of type $\D_n$ is shown to
be semisimple and free of rank $(2^n+1)n!!-(2^{n-1}+1)n!$ over a specified
commutative ring $R$, where $n!! = 1\cdot 3\cdots (2n-1)$.  We also show it
is a cellular algebra over suitable ring extensions of $R$.  The Brauer
algebra of type $\D_n$ is the image of an $R$-equivariant homomorphism
and is also semisimple and free of the same rank, but over the ring
$\Z[\delta^{\pm1}]$. A rewrite system for the Brauer algebra is used in
bounding the rank of the BMW algebra above.  As a consequence of our
results, the generalized Temperley--Lieb algebra of type $\D_n$ is a
subalgebra of the BMW algebra of the same type.
\end{abstract}

\maketitle

\medskip
{\sc keywords:} associative algebra, Birman--Murakami--Wenzl algebra, BMW
algebra, Brauer algebra, cellular algebra, Coxeter group, generalized
Temperley--Lieb algebra, root system, semisimple algebra, word problem in
semigroups


\section{Introduction}
\label{sec:introduction}
In \cite{BirWen}, Birman and Wenzl, and independently in \cite{Mur},
Murakami, defined algebras indexed by the natural numbers which play a role
in both the representation theory of quantum groups and knot theory.  They
were given by generators and relations.  In \cite{MorWas}, Morton and
Wasserman gave them a description in terms of tangles.  These are the
Birman--Murakami--Wenzl algebras (usually abbreviated to BMW algebras) for
the Coxeter system of type $\A_{n}$.  They behave nicely with respect to
restriction to the algebras generated by subsets of the generators.  For
instance, the BMW algebras of a restricted type embed naturally into the
bigger ones.  This is similar to the fact that in Weyl groups subgroups
generated by subsets of the standard reflections are themselves Weyl groups.
The Hecke algebra of type $\A_{n}$ is a natural quotient of the
Birman--Murakami--Wenzl algebra of type $\A_{n}$ and the Temperley--Lieb
algebra, conceived originally for statistics (cf.~\cite{Temp-Lieb}), is a
natural subalgebra.  Inspired by the beauty of these results, the existence
of Temperley--Lieb algebras of other types (\cite{Fan,Gra,GrLe,Gre}) and the
existence of a faithful linear representation of the braid group
(\cite{CohWal,digne}), the authors defined analogues for other simply laced
Coxeter diagrams and found some of their properties in \cite{CGW}.  The
faithful linear representations of the braid group were shown first by
Bigelow in \cite{Bigelow} and Krammer in \cite{Krammer}.  They used a
representation introduced by Lawrence in \cite{Lawrence}.

In this paper we consider the algebras when the Coxeter diagram is of type
$\D_n$. We prove the conjecture stated in \cite[Section 7.1]{CGW}, which is
Theorem \ref{th:main}.  Here, $n!! = 1\cdot 3\cdots (2n-1)$.  We work over
the quotient ring $R$ of $\Z[\delta,\delta^{-1},l,l^{-1},m]$ by the ideal
generated by $m(1-\delta)-(l-l^{-1})$ instead of the field $\Q(l,\delta)$ in
which it embeds (see Lemma~\ref{lm:RinQldelta}).

\begin{Thm}\label{th:main}
The BMW algebra of type $\D_n$ over $R$ is free of rank $$(2^n+1)
n!!-(2^{n-1}+1)n!.$$ When tensored with $\Q(l,\delta)$, it is semisimple.
\end{Thm}

\np The result produces linear representations of the Artin group of type
$\D_n$ similar to the representations of the braid group on $n$ strands which
arose from the BMW algebra of type $\A_{n-1}$.  These include the faithful
representations related to the Lawrence--Krammer representations occurring in
\cite{CohWal} as well as the representations occurring in \cite{CGW}.
Furthermore, specific information about the representations is given in terms
of sets of orthogonal roots and irreducible representations of Weyl groups of
type $\D_r$ for certain $r$ (cf.~Remark \ref{rmk:reps}).

These sets of orthogonal roots are also used in the description of a
cellular basis, whose elements are determined by pairs of such root sets and
a Weyl group element.
This leads, for suitable extensions of the coefficient ring $R$,
to cellularity of the BMW algebra
$\BMW(\D_n)$ in the sense of \cite[Definition~1.1]{GL}. For $\BMW(\A_n)$,
this result is known thanks to \cite{Xi}.

\begin{Thm}\label{th:cellular}
The BMW algebra of type $\D_n$ is cellular if the coefficient
ring $R$ is extended to an integral domain containing
an inverse to $2$.
\end{Thm}

As a consequence of the work we are able to show the Temperley--Lieb algebra
of type $\D_n$ as defined in \cite{Fan,Gra,Gre}  is a natural subalgebra.

\begin{Cor}\label{TempLieb}  The generalized Temperley--Lieb algebra of type $\D_n$ is a natural
subalgebra of both the Brauer algebra and the BMW algebra of type $\D_n$
over the rings $Z[\delta,\delta^{-1}]$ and $R$, respectively.
\end{Cor}

The current work completes the proof that there is an isomorphism from the
BMW algebra to the algebra of tangles having a pole of order 2 studied in
\cite{CGW3}. For each element of the cellular basis, the two corresponding
root sets determine the set of horizontal strands at the top and bottom,
respectively, and the corresponding Weyl group element determines the
vertical strands of the tangle. The isomorphism is discussed at the end of
this paper.

Putting together \cite{MorWas}, Theorem \ref{th:main} and the main theorem of \cite{CW4},
we have reached a complete description of the BMW algebras of spherical
simply laced type.

\section{Overview}
\label{sec:overview}
We proceed as follows. First, in Section \ref{BMWrelations}, we introduce
the BMW algebra $\BMW(\M)$ over $R$ for $\M$ of type $\A_n$ $(n\ge1)$,
$\D_n$ $(n\ge4)$, or $\E_n$ $(n=6,7,8)$, which we denote $\ADE$.  Then the
Brauer algebra, $\Br(\M)$, of the same type over $\Z[\delta^{\pm1}]$ is
obtained from $\BMW(M)$ by specializing $m$ to $0$ and $l$ to $1$.  This
algebra was defined in \cite{CFW} where it was shown to be free over $R$ of
rank $(2^n+1) n!!-(2^{n-1}+1)n!$ in case $M=\D_n$. The modding out of $m$
and $l-1$ gives a surjective $R$-equivariant map $\mu:\BMW(\M) \mapsto
\Br(\M)$.

The Brauer algebra $\Br(M)$ is given in terms of generators $e_i$, $r_i$ for
$i$ running over the nodes of $M$, and relations determined by $M$
(cf.~Definition \ref{df:BrMonoid}). The
subalgebra of $\Br(M)$ generated by the $r_i$ is the group algebra over
$\Z[\delta^{\pm1}]$ of $W(M)$, the Coxeter group of type $M$.

The specialization enables us to pass from monomials in $\BMW(\D_n)$ to
monomials in $\Br(\D_n)$.  We will use this observation to find a basis of
monomials for $\BMW(\D_n)$ from a similar basis in $\Br(\D_n)$.

In Section~\ref{sec:functionMonoiod} we summarize results from \cite{CFW}
and \cite{CGW2} which show how the monomials of $\Br(M)$ determine sets of
mutually orthogonal roots, which in the case $M=\A_{n-1}$ are directly
related to tops and bottoms of the well-known Brauer diagrams. The
monomials, including powers of $\delta$, form a monoid inside $\Br(M)$,
denoted $\BrM(M)$ (see Definition \ref{df:BrMonoid}).

In Sections~\ref{sec:reduction} and \ref{sec:reduction2} we use the
following strategy to produce a basis of $\BMW(\D_n)$ from elements of
$\BrM(\D_n)$.  A word $\und a$ in the generators of the Brauer monoid
$\BrM(M)$ is said to be of \emph{height} $t$ if the number of generators
$r_i$ occurring in it is equal to $t$.  We say that $\und a$ is
\emph{reducible} to another word $\und b$ if $\und b$ can be obtained from
$\und a$ by a finite sequence of specified rewrite rules (listed in Table
\ref{BrauerTable}) that do not increase the height.  This process will be
called a \emph{reduction}.  The significance of such a reduction is that the
word $\und a$ also corresponds to a unique monomial in the BMW algebra and
that a parallel reduction (with rules listed in Table \ref{BMWTable}) can be
carried out in the BMW algebra in the sense that the monomial in
$\BMW(\D_n)$ corresponding to $\und a$ can be rewritten as a linear
combination of monomials all of which are represented by words of height
less than or equal to the height of $\und a$, with equality occurring for at
most one term (see Proposition \ref{prop:wordssameinbrauer}(ii)).  We
exhibit a finite set of reduced words to which each word reduces; see
Corollary~\ref{cor:finalreduction}.  This will lead to a set $T$ of reduced
words such that every word in the generators of $\BMW(\D_n)$ can be reduced
to an element of $T$ up to multiples by powers of $\delta$.  The above
argument will give that, when viewed as elements of $\BMW(\D_n)$, the set
$T$ is a spanning set of $\BMW(\D_n)$.

In Section~\ref{sec:conclusion} we prove our main result by constructing a
suitable set $T$ of monomials corresponding to specific triples consisting
of pairs of sets of mutually orthogonal roots and a Weyl group element. We
also prove Corollary~\ref{TempLieb} by showing that the generalized
Temperley--Lieb algebra of type $\D_n$, embeds in $\BMW(\D_n)$ and in
$\Br(\D_n)$.

In Section~\ref{cellular} we show that if the ring of coefficients is
extended to an integral domain containing $2^{-1}$, the algebra $\BMW(\D_n)$
is cellular in the sense of \cite[Definition~1.1]{GL}.  In our proof, we
need the ring extension in order to invoke \cite[Theorem~1.1]{MG} where
cellularity of the Hecke algebras of type $\D_n$ is proved for such rings of
coefficients.  This Hecke algebra is a natural quotient of $\BMW(\D_n)$ and
the Hecke algebras of type $\D_{n-2t}$ occur as subalgebras with different
idempotents as identities in the analysis.

We have applied the above results in \cite{CGW3}, where a tangle algebra
$\KT(\D_n)$ over $R$ on $n$ strands was introduced. This algebra was shown
to be a homomorphic image of the BMW algebra $\BMW(\D_n)$ of type $\D_n$ and
Theorem \ref{th:main} gives that $\KT(\D_n)$ is an isomorphic image of it.

Part of the work reported here grew out of the PhD.~thesis of one of us,
\cite{DAHG}. The other two authors wish to acknowledge Caltech and Technische
Universiteit Eindhoven for enabling mutual visits.

\section{BMW and Brauer Algebras}\label{BMWrelations}
The BMW algebras of type $\A_n$ $(n\ge1)$, $\D_n$ $(n\ge4)$, and $\E_n$
$(n=6,7,8)$ have been discussed extensively in \cite{CGW}.  We assume that
$\M$ is a Coxeter diagram which is one of these (in particular, it has no
multiple bonds).  Our main results will only concern $\M$ of type $\A_{n-1}$
and $\D_n$.  The BMW algebra of type $\M$ is defined over the ring $R =
\Z[l^{\pm1},m,\delta^{\pm1}]/(m(\delta-1)-(l^{-1}-l))$.

\begin{Def}\label{df:BMW}\rm
The BMW algebra $\BMW(\M)$ of type $\M$ is the free algebra over $R$ given by
generators $g_i$, $e_i$ with $i$ running over the nodes of the diagram $M$,
subject to the relations in the BMW Relations Table \ref{BMWTable} where
$i \sim j$ denotes adjacency of two
nodes $i$ and $j$.

\begin{center}
\begin{table}[ht]
\begin{tabular}{|lcl|}
\hline
&& for $i$\\
(RSrr)&\quad&$g_i^2=1-m(g_i-l^{-1}e_i)$\\
(RSer)&\quad&$e_ig_i=l^{-1}e_i$\\
(RSre)&\quad&$g_ie_i=l^{-1}e_i$\\
(HSee)&\quad&$e_i^2=\delta e_i$\\
\hline
&&for $i\not\sim j$\\
(HCrr)&\quad&$g_ig_j=g_jg_i$\\
(HCer)&\quad&$e_ig_j= g_je_i$\\
(HCee)&\quad&$e_ie_j= e_je_i$\\
\hline
&&for $i\sim j$\\
(HNrrr)&\quad&$g_ig_jg_i=g_jg_ig_j$\\
(HNrer)&\quad&$g_je_ig_j=g_ie_jg_i+ m(e_jg_i - e_ig_j + g_ie_j - g_je_i)$\\
&&$\qquad\qquad + m^2(e_j-e_i)$\\
(RNrre)&\quad&$g_jg_ie_j=e_ie_j$\\
(RNerr)&\quad&$e_ig_jg_i=e_ie_j$\\
(HNree)&\quad&$g_je_ie_j=g_ie_j+m(e_j-e_ie_j)$\\
(RNere)&\quad&$e_ig_je_i=le_i$\\
(HNeer)&\quad&$e_je_ig_j=e_jg_i+m(e_j-e_je_i)$\\
(HNeee)&\quad&$e_ie_je_i=e_i$\\
\hline
&&for $i\sim j\sim k$\\
(HTeere)&\quad&$e_je_ig_ke_j=e_jg_ie_ke_j$\\
(RTerre)&\quad&$e_jg_ig_ke_j=e_je_ie_ke_j + m (e_je_ig_ke_j - l e_j)$\\
\hline
\end{tabular}
\bigskip
\caption{BMW Relations}
\label{BMWTable}
\end{table}
\end{center}
\end{Def}

\begin{Remark}
\rm The set of relations given is redundant. In fact, the relations (HNrer),
(HNree), (HNeer), (HNeee), (HTeere), and (RTerre) follow from the others, as
we will explain.  Moreover, if $\BMW(M)$ is tensored with a ring in which
$m$ is invertible, they all follow from (RSrr), (RSre), (HCrr), (HNrrr), and
(RNere).  This is shown in \cite{CGW} where these were labeled (D1), (R1),
(B1), (B2), and (R2), respectively.

We will prove the stated redundancies, starting with (HNeee).  By (RNerr),
(RSre), and (RNere), respectively, we have
\begin{eqnarray*}
e_ie_je_i&=&e_ig_jg_ie_i = l^{-1}e_ig_je_i   =e_i.
 \end{eqnarray*}
For (HNeer) we multiply (RNerr) from the right by $g_j$, apply (RSrr) to the
right hand side, and,
for the final equality, (RNerr), and (RNere):
\begin{eqnarray*}
e_je_ig_j  &=&  e_jg_ig_j^2
           =e_jg_i(1-mg_j+ml^{-1}e_j)
            =e_jg_i-me_jg_ig_j+ml^{-1}e_jg_ie_j  \\
            &=& e_jg_i-me_je_i+me_j.
\end{eqnarray*}
(HNree) is derived in a similar way.
The equation (HNrer) is dealt with in \cite[Proposition 2.3]{CGW}
by use of the relations we have obtained.
For (HTeere), we use (RNerr) and (RNrre), respectively:
$$e_je_ig_ke_j=e_jg_ig_jg_ke_j=e_jg_ie_ke_j.$$  Recall here that
$i \not \sim k$ because the diagram $M$ has no triangles.
For (RTerre) write $e_jg_ig_ke_j=e_jg_ig_jg_j^{-1}g_ke_j$
and use the expression $g_j^{-1}=g_j+m-me_j $ which follows
from (RSrr) and is given in \cite[Proposition 2.1]{CGW}.
\end{Remark}

\begin{Defs}\label{df:BrMonoid}\rm
Let $\M$ be a graph of type $\ADE$.
We define the Brauer monoid $\BrM(\M)$ to be the monoid
generated by the elements $r_i$ and $e_i$ $(i\in \M)$ and $\delta$
subject to the relations in the Brauer Relations Table~\ref{BrauerTable}.
The Brauer algebra of type $\M$ is the monoid algebra $\Z[\BrM(\M)]$.

\begin{center}
\begin{table}[ht]
\begin{tabular}{|lcl|lcl|}
\hline
label&\quad&relation&label&\quad&relation\\
\hline
($\delta$)&\quad&$\delta$ is central&($\delta^{-1}$)&\quad&$\delta\delta^{-1}=1$\\
\hline
&\multispan{4}{\hfill for $i$\hfill}&\\
\hline
(RSrr)&\quad&$r_i^2=1$&
(RSer)&\quad&$e_ir_i=e_i$\\
(RSre)&\quad&$r_ie_i=e_i$&
(HSee)&\quad&$e_i^2=\delta e_i$\\
\hline
&\multispan{4}{\hfill for $i\not\sim j$\hfill}&\\
\hline
(HCrr)&\quad&$r_ir_j=r_jr_i$&
(HCer)&\quad&$e_ir_j= r_je_i$\\
(HCee)&\quad&$e_ie_j= e_je_i$&&&\\
\hline
&\multispan{4}{\hfill for $i\sim j$\hfill}&\\
\hline
(HNrrr)&\quad&$r_ir_jr_i=r_jr_ir_j$&
(HNrer)&\quad&$r_je_ir_j=r_ie_jr_i$\\
(RNrre)&\quad&$r_jr_ie_j=e_ie_j$&
(RNerr)&\quad&$e_ir_jr_i=e_ie_j$\\
(HNree)&\quad&$r_je_ie_j=r_ie_j$&
(RNere)&\quad&$e_ir_je_i=e_i$\\
(HNeer)&\quad&$e_je_ir_j=e_jr_i$&
(HNeee)&\quad&$e_ie_je_i=e_i$\\
\hline
&\multispan{4}{\hfill for $i\sim j\sim k$\hfill}&\\
\hline
(HTeere)&\quad&$e_je_ir_ke_j=e_jr_ie_ke_j$&
(RTerre)&\quad&$e_jr_ir_ke_j=e_je_ie_ke_j$\\
\hline
\end{tabular}
\bigskip
\caption{Brauer Relations}
\label{BrauerTable}
\end{table}
\end{center}

The $r_i$ in $\BrM(\M)$ generate a subgroup of the Brauer monoid that we
denote $W$.  This is a Coxeter group of type $\M$ as the $r_i$ satisfy the
required relations and, after factoring out the ideal of $\Br(\M)$ generated
by the $e_i$, we obtain the group algebra of $W$ over $\Z[\delta^{\pm1}]$.

We consider the Brauer algebra of type $\M$ as an algebra over
$\Z[\delta^{\pm1}]$.  Here $\delta$ is in the center of $\BrM(\M)$ and we
identify this $\delta $ with the $\delta $ in $Z[\delta^{\pm1}]$.  Since the
other defining relations of the Brauer monoid are the defining relations of
the corresponding BMW algebra $\BMW(\M)$ modulo the ideal $(l-1,m)$
generated by $l-1$ and $m$, the Brauer algebra $\Br(M)$ can be identified
with $\BMW(\M)\otimes_R R/(l-1,m)$. The corresponding equivariant map
$a\mapsto a\otimes 1$ will be denoted by $\mu$.
\end{Defs}

Just as for $\BMW(M)$, some of the relations in the Brauer Relations Table
\ref{BrauerTable} are redundant; see \cite[Lemma 3.1]{CFW}.  We will need to
rewrite words in the generators $r_i$ and $e_i$, with $\delta^{\pm1}$
viewed as coefficients. This necessitates the extra relations that are
displayed in Table \ref{BrauerTable}.

\begin{Defs}\label{df:reduction} \rm
By $F_n$ we denote the monoid that is the central product of the free monoid
on the symbols $r_i$, $e_i$ $(i=1,\ldots,n)$ with the infinite cyclic group
generated by $\delta$. Its elements will be called \emph{words}.  There is a
surjective homomorphism of monoids $\pi : F_n\to \BrM(\M)$ mapping the
symbols $r_i$, $e_i$, and $\delta$ to the corresponding elements of
$\BrM(\M)$. The monomial in $\BMW(\M)$ corresponding to $\und a\in F_n$,
obtained by replacing $r_i$ by $g_i$ and leaving $e_i$ and $\delta$ as
before, will be denoted $\rho(\und a)$, so $\mu(\rho(\und a)) = \pi(\und
a)$. A word $\und a\in F_n$ is said to be
of {\em height} $t$ if the number of $r_i$ occurring in it is equal to $t$;
we denote this number $t$ by $\het(\und a)$.

We say that $\und a$ is \emph{reducible} to another word $ \und b$, or that
$\und b$ is a \emph{reduction} of $\und a$, if $ \und b$ can be obtained by
a finite sequence of specified rewrites, listed in the Brauer Relations
Table~\ref{BrauerTable}, starting from $\und a$, that do not increase the
height.  We call a word in $F_n$ \emph{reduced} if it cannot be further
reduced to a word of smaller height.  We have labeled the relations in the
tables above with R or H according to whether the rewrite from left to right
strictly lowers the height or not.  If the number stays the same, we call it
H for \emph{homogeneous}. Our rewrite system will be the set of all rewrites
in the Brauer Relations Table \ref{BrauerTable} in either direction in the
homogeneous case when an H appears in its label and from left to right only
in case an R occurs in its label.  We write $\und a\homog \und b$ if $\und
a$ can be reduced to $\und b$; for example (RNerr) gives $e_2r_3e_2\homog
e_2$ if $2\sim 3$.  If the height does not decrease during a reduction, we
sometimes use the term {\em homogeneous reduction} and write $\und a \isog
\und b$; for example, (HNeee) gives $e_2\isog e_2e_3e_2$ if $2\sim 3$. If it
does decrease, we also speak of a {\em strict reduction}.

Homogeneous reduction induces a congruence relation on $F_n$, to which we
will refer as {\em homogeneous equivalence}.  We denote its set of
equivalence classes by $F_n/\isog$. The congruence property turns it into a
monoid.
\end{Defs}

\np
The reductions in $F_n$ are important
because they have a meaning for both the Brauer algebra and the corresponding
BMW algebra.  For each of the relations in the Brauer Relations Table
\ref{BrauerTable}, there is a corresponding relation in the BMW Relations
Table \ref{BMWTable}.  In Section~\ref{sec:conclusion}, the following
proposition will be used to find a basis of $\BMW(\D_n)$ that has the same
size as a basis of $\Br(\D_n) $.

\suppress{**amc: this has been said in Definitions 3.3 already, therefore
removed here:
\begin{Remark}  \rm When working with the Brauer monoid we  identify the $\delta$ in the Brauer monoid with the $\delta$
in $Z[\delta^{\pm 1}]$ and when working in the BMW algebra we identify the $\delta$ in the BMW algebra with the $\delta $ in $R$.
\end{Remark}
}

\begin{Prop}\label{prop:wordssameinbrauer}
Suppose $\und a\homog \und b$ with $\und a,\ \und b \in F_n$.
\begin{enumerate}[(i)]
\item
$\pi(\und a)=\pi (\und b)$ in
$\BrM(\M)$.
\item\label{prop:BMWBrRewrite} There are a finite number of $\lambda_{\und
c}\in R$ such that, in $\BMW(M)$, $$\rho(\und a) = \rho(\und b) + \sum_{\und
c\in F_n,\,\het(\und c)<\het(\und a)}m\lambda_{\und c}\rho(\und c).$$
\end{enumerate}
\end{Prop}

\nl
\begin{proof}
(i).  For each reduction step of the sequence of relations, the word
evaluated in $\BrM(\M)$ is the same because the relations are satisfied in
$\BrM(\M)$ by definition.  This means $\pi(\und a)=\pi(\und b)$, proving
(i).

\nl(ii).  The expressions in the BMW Relations Table \ref{BMWTable} all have
one term on each side whose coefficient is not a multiple of $m$.  These
terms are the same as in the Brauer Relations Table \ref{BrauerTable} with
$g_i$ instead of $r_i$.  Indeed, if $l=1$ and the terms with coefficient $m$
are ignored, the tables are the same.  Each reduction step in $\und a\homog
\und b$, replaces the term on the left with the corresponding one on the
right side of the equality in the table plus terms that are multiples of $m$
and have strictly smaller height.  The end result is $\rho (\und b) $ plus
terms that are multiples of $m$, whose height has been reduced at least
once.  As $\underline{a}\homog\underline{b}$ involves only a finite sequence
of specific rewrites from Table \ref{BrauerTable}, only a finite number of
substitutions from Table \ref{BMWTable} has been applied, and so only a
finite number of summands occurs at the right hand side of the equality
in (ii).
\end{proof}

\begin{Notation}
\label{not:op}
\rm For $x_1,\ldots, x_q\in\{r_1,\ldots,r_n,e_1,\ldots,e_n,\delta^{\pm1}\}$,
we write $(x_1\cdots x_q)^\op = x_q\cdots x_1$, thus defining an opposition
map on $F_n$. This notation is compatible with the maps $\pi$ and $\rho$
when ${\cdot}^\op$ on $\BMW(\M)$ and $\Br(\M)$ is interpreted as the
anti-involution of \cite[Remark~2.1(i)]{CGW} and \cite[Remark~5.7]{CFW},
respectively.
\end{Notation}

\np
To end this section, we discuss properties of $R$ which show how to relate
some properties of sets of monomials in $\BMW(\D_n)$ to corresponding ones
in $\Br(\D_n)$ using the maps $\pi$ and $\rho$.

\begin{Lm}\label{lm:RinQldelta}
The ring $R$ embeds in $\Q(\delta)[l^{\pm1}]$ and also in $\Q(l,\delta)$.
\end{Lm}

\nl
\begin{proof}
Let $D=\Z[l^{\pm1},\delta^{\pm1}]$, which is a unique factorization domain,
and let $F$ be its field of fractions.  Put $s(m)=(1-\delta)m-(l-l^{-1})$.
Notice $s(m) $ is primitive and so irreducible in $F[m]$ by Gauss' Lemma.
Hence $R = D[m]/(s(m))$ is an
integral domain. Its field of fractions is $\Q(l,\delta)$.
Finally, $\Q(\delta)[l^{\pm1}]$ is a subring of
$\Q(l,\delta)$ containing both $D$ and $m$, as the latter
is equal to $(l-l^{-1})/(1-\delta)$ modulo
$s(m)$, and so also contains $R$.
\end{proof}

\np
The following lemma will give a lower bound for the rank of $\BMW(\D_n)$.

\begin{Lm}\label{lm:basis}  Suppose that
$T$ is a finite set of monomials in $F_n$ whose images
$(\pi(t))_{t\in T}$ are
linearly independent in $\Br(\D_n)$.
Then $(\rho(t))_{t\in T}$ are linearly independent in $\BMW(\D_n)$.
\end{Lm}

\nl
\begin{proof}
Suppose that $\sum_{t\in T}\lambda_t\rho(t)$ with $\lambda_t\in R$ is a
non-trivial linear combination that is equal to $0$ in $\BMW(\D_n)$.  Then
the same non-trivial linear relation holds over the principal ideal domain
$\Q(\delta)[l^{\pm1}]$ into which $R$ embeds according to
Lemma~\ref{lm:RinQldelta}.  Rescale the coefficients by a suitable power of
$l-1$ to guarantee $\lambda_s\not\in (l-1) \Q(\delta)[l^{\pm1}]$ for some
$s\in T$.  Now $\mu(\lambda_s)\ne 0$ and $\pi(t) = \mu(\rho(t))$ for $t\in T$
(cf.~Definitions \ref{df:reduction}),
so $\sum_{t\in T}\mu(\lambda_t)\pi(t)$ is a non-trivial linear combination in
$\Br(\D_n)$, that is equal to $0$, contradicting the linear independence
assumption on $\pi(t)_{t\in T}$.
\end{proof}

\np
The following result will yield the right upper bound on the rank of
$\BMW(\D_n)$.

\begin{Prop}\label{prop:BMWbasis}
Let $M$ be of type $\ADE$ and let $T$ be a set of words in $F_n$ such that
the $(\pi(t))_{t\in T}$ is a basis of $\Br(M)$.  If each
word in $F_n$ can be reduced to an element of $\delta^{\Z}T$, then $\rho(T)$
is a basis of $\BMW(M)$ and each element of $T$ is a reduced word.
\end{Prop}

\nl
\begin{proof}
Assume that each word in $F_n$ can be reduced to an element of
$\delta^{\Z}T$.  We first prove that $\rho(T)$ is a linear spanning set of
$\BMW(M)$. If not, there is a word $\und a$ in $F_n$ such that $\rho(\und
a)$ is not in the linear span of $\rho(T)$. Pick one of smallest
height. Then, by assumption, $\und a\homog \und b$ for some $\und b\in
\delta^{\Z}T$.
Proposition~\ref{prop:wordssameinbrauer}(\ref{prop:BMWBrRewrite}) implies
that $\rho(\und a)-\rho(\und b)$ is a linear combination of monomials in
$\BMW(M)$ of height lower than $s=\het(\und a)$.  If $s = 0$, this means
$\rho(\und a) = \rho(\und b)\in\delta^{\Z}\rho(T)$.  Otherwise $s>0$ and we may assume,
using induction on height, that monomials in $\BMW(M)$ of height lower than
$s$ are all in the linear span of the elements in $\rho(T)$ of height lower
than $s$.  Then the right hand side in the expression of $\rho(\und
a)-\rho(\und b)$ as a linear combination of monomials of lower height is in
the linear span of $\rho(T)$.  Consequently, $\rho(\und a)$ is in the same
linear span, a contradiction.  We have shown that $\BMW(M)$ is spanned by
$\rho(T)$.

It now follows from Lemma \ref{lm:basis} that $(\rho(t))_{t\in T}$ is
a basis of $\BMW(\D_n)$. If $\und{a}\in T$ is not reduced, then there
is $\und{b}\in F_n$ with $\und a\homog \und b$ and $\het(\und a) >
\het(\und b)$. After applying the assumption to $\und{b}$, we may
assume $\und b \in\delta^{\Z}T$ and still $\het(\und a)>\het(\und b)$.
In view of the hypothesis, $\pi(\und a)$ and $\delta^i\pi(\und b)$
for some $i\in \Z$ are
distinct members of a basis of $\Br(\D_n)$,
so $\pi(\und a)$ and $\pi(\und b)$
are linearly independent.  On
the other hand, by Proposition~\ref{prop:wordssameinbrauer}(i),
$\pi(\und a) = \pi(\und b)$, a contradiction.
\end{proof}

\section{Admissible Sets and the Function Monoid}\label{sec:functionMonoiod}
\begin{center}
\begin{figure}[ht]
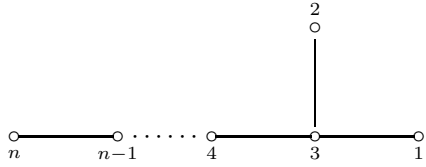

$$\Dnrev$$
\caption{The diagram of type $\D_n$ with node labels}
\label{fig:Dn}
\end{figure}
\end{center}

Let $n\in\N$, $n\ge4$. We summarize some of the results of \cite{CFW}
and \cite{CGW2} about admissible sets with a special focus on type $M =
\D_n$.  These are particular sets of mutually orthogonal positive roots.
The results will be used to monitor the reduction of words in $F_n$.  We
will fix a root system $\Phi$ for $W$ and a set of simple roots
$\a_1,\ldots,\a_n$ with indices for $M = \D_n$ as indicated in the Dynkin
diagram of Figure \ref{fig:Dn}.

\noindent
In terms of the standard orthonormal basis $\eps_1,\ldots,\eps_n$ of $\R^n$,
these simple roots are $\a_1=\eps_1+\eps_2$, $\a_2=\eps_2-\eps_1,$
$\a_3=\eps_3-\eps_2,$ $\ldots$, $\a_n=\eps_n-\eps_{n-1}$.  Accordingly, we
will write $\Phi^+ = (\Z_{\ge0}\a_1+\Z_{\ge0}\a_2+\cdots+\Z_{\ge0}\a_n)\cap
\Phi$ where $\Z_{\ge0}$ are the non-negative integers.
The elements of $\Phi^+$ are called the \emph{positive roots} of $\Phi$ (or
simply $\D_n$); they are of the form $\eps_j-\eps_i$ and $\eps_i+\eps_j$ for
$n\geq j > i\geq 1$.  Recall $\Phi=\Phi^+\cup (-\Phi^+)$.  The reflection in
$\R^n$ with root $\beta$ is denoted $r_\beta$. The map $r_i\mapsto
r_{\alpha_i}$ $(i=1,\ldots,n)$ extends to an isomorphism from $W$ to a
reflection subgroup of the orthogonal group on $\R^n$.  We often identify
$W$ with this reflection group by means of the isomorphism.

There are some standard properties of the root systems we are using which we
mention here for convenience.  All roots have square norm $2$.  The inner
products are all $\pm 2$, $\pm1$, or $ 0$.  If $(\b,\c)=1$, then $\b-\c$ is
a root, and if $(\b,\c)=-1$, then $\b+\c$ is a root, equal to $r_\b\c$.
Further if $(\b,\c)=0$, then $\b\pm \c$ is never a root.  We often encounter
the situation in which $(\b,\a_i)=0$, $i\sim j$, and $\b-\a_j$ is a root.
Then $(\b-\a_j, \a_i)=1$ and so $\b-\a_j-\a_i = r_ir_j\b$ is also a
root.

\begin{Remark}
\label{rmk:height}
\rm
There are two notions of height.  The first is $\het (\und a)$ for $\und a$
an element of $F_n$ (cf.~Definitions \ref{df:reduction}). The second is the
more standard notion of height of a positive root $\b$.  This is $\sum \lambda_i$
for the root $\sum \lambda_i \a_i$ where the $\a_i$ are the simple roots.  We also
denote this $\het (\b)$ and trust no confusion will arise.
\end{Remark}

In order to recognize the elements of the ideal in $\BrM(\D_n)$
generated by $e_1e_2$, see Definition~\ref{df:Theta}, we will need the
notion of orthogonal mates.

\begin{Defs}\label{def:rn*}
\rm For $\b=\eps_i-\eps_j$ a root in the root system $\Phi$ of type $\D_n$
embedded in $\R^n$ as indicated above, its \emph{orthogonal mate} is defined
to be $\beta^* = \eps_i+\eps_j$ and, vice versa, the \emph{orthogonal mate}
of $\beta^*$ is $\beta^{**} = \beta$.  Furthermore, we write $r_\b^*$ for
$r_{\b^*}$, the reflection whose root is the orthogonal mate of $\b$.  For
the simple roots $\a_i$ we also write $r_i^*$ instead of $r_{\a_i}^*$.

If $n>4$, the roots orthogonal to $\beta$ form a subsystem of
$\Phi$ of type $\A_1\D_{n-2}$ and $\beta^*$ is the unique positive root in
the $\A_1$ component of this subsystem.  If $n=4$, the choice of orthogonal
mate essentially depends on the choice of an orthogonal pair of simple
roots.

There are several equivalent definitions of admissible sets as outlined
in \cite[Proposition 2.3]{CGW2}.  For our purposes we may define a set $B$
of mutually orthogonal positive roots to be {\em admissible} if and only if,
when $\a_1,\a_2,\a_3\in B$ and there exists a root $\a$ for which
$(\a_i,\a)=\pm1$ for all $i$, then $r_{\a}r_{\a_1}r_{\a_2}r_{\a_3}\a$ or
$-r_{\a}r_{\a_1}r_{\a_2}r_{\a_3}\a$ is also in $B$.  Given any set, $B$, of
mutually orthogonal positive roots, a straightforward exercise shows there
is a unique smallest admissible set containing $B$.  This set is called the
{\em admissible closure} of $B$, notation $\cl{B}$; see
\cite[Definition~2.2]{CFW}.

\label{df:heightAO}
By $\AO$ we denote the collection of all admissible sets (including the
empty set). This set has a natural $W$-action given by
$$wB = \Phi^+\cap\left\{ \pm w\b \mid\b\in B \right\}$$ for $w\in W$.  A
representative of each $W$-orbit in $\AO$ is given in \cite[Table 3]{CFW};
this is a corrected version of a similar table in \cite{CGW2}.  We will need
these only for types $\A_n$ and $\D_n$, which, for the convenience of the
reader, are summarized in Lemma \ref{lm:Worbits} and Table
\ref{table:types2}.
The meaning of $M_Y$ and $S_Y$ in Table \ref{table:types2} will become clear
later (in Proposition \ref{prop:ZY}).
\end{Defs}

\begin{Notation}\label{def:theta}  \rm
By $\MY$ we denote the collection of the following sets of nodes of $\D_n$.
\begin{eqnarray*}
Y(t) &=& \{ n,n-2,\ldots,n-2t+2\}, \mbox{   for   } t\in
[0,\lfloor n/2\rfloor]\\
Y^*(t) &=& \{n,n-2,\ldots,n-2t+4, 1, 2\}, \mbox{   for   }
t\in[1,\lfloor n/2\rfloor]\\
Y'(n/2) &=& \{ 1,4,6,8,\ldots,n\} \mbox{ if  } n \mbox{ is even}.
\end{eqnarray*}

For $Y\in\MY$, the set $B_Y$ is the admissible closure of the set
of roots $\a_j$ for $j\in Y$.
For $Y=Y^*(t)$,
this implies that $B_Y$ is the set of roots $\a_j$ and $\a_j^*$
for $j\in Y$.
For $Y=Y(t)$ or $Y'(n/2)$ however,
no orthogonal mates occur
and so $B_Y$ is the set of roots $\a_j$
for $j\in Y$.
\end{Notation}

\begin{Lm}
\label{lm:Worbits}
Each $W$-orbit in $\AO$ has a unique representative $B_Y$ for $Y\in\MY$.
\end{Lm}

\np For instance, if $n=4$, there are three orbits of admissible sets of
size $2$, with representatives $B_Y$, where $Y = Y(2)$, $ Y'(2)$, and
$Y^*(1)$, respectively.

Notice that, if $B$ is an admissible set
containing a root as well as its orthogonal mate, then it is a union of
roots together with their orthogonal mates.

\bigskip
\begin{table}[ht]
\begin{center}
\begin{tabular}{|c|cccc|}
\hline
${\CoxDiag}$ & $|B_Y|$ &  $M_Y$ & $Y$&$S_Y=\{x\he_Y\mid
x\mbox{ as below}\}$\\
\hline
$\A_n$ & $t$ & $\A_{n - 2t}$ & $Y(t)$ & $r_1,\ldots,r_{n-2t}$\\
\hline
$\D_n$ & $t$&
$\A_1 \D_{n - 2t}$ & $Y(t)$&$\und z_n^*, r_1,\ldots,r_{n-2t}$\\
$\D_n$& $0$  & $ \D_{n}$ & $\emptyset$&$r_1,\ldots,r_n$\\
$\D_n$ ($n$ odd)&$(n-1)/2$  & $\A_1$ & $Y((n-1)/2)$&$\und z_n^*$\\
$\D_n$ ($n$ even)&$n/2$ & $\A_1$ & $Y(n/2)$&$r_1$\\
$\D_n$ ($n$ even)& $n/2$ & $\A_1$ & $Y'(n/2)$&$r_2$\\
$\D_n$ & $2t$ & $\A_{n - 2t - 1}$ & $Y^*(t)$ &
$r_4,\ldots,r_{n-2t+2}$\\
\hline
\end{tabular}
\end{center}
\bigskip
\caption{\label{table:types2}\textrm{Cocliques $Y\in\MY$ of
$M$ and sizes of the admissible sets $B_Y$ for types $\A_n$ and $\D_n$.
The index $t$ in the line for $\A_n$ satisfies
$0\le t\le \lfloor n/2\rfloor$.
The index $t$ in the first line for $\D_n$ satisfies
$0< t< \lfloor (n-1)/2\rfloor$. In the last line it satisfies
$0< t\le \lfloor n/2\rfloor$.
}}
\end{table}

The following proposition is proved
in \cite[Theorem 3.6]{CFW}; the fact that $e_iB$ as described below is well
defined is shown in \cite[Lemma 3.3(v)]{CFW}.

\begin{Prop}\label{prop:BrMonAction}
Let $M$ be of type $\ADE$. The action of $W$ on $\AO$ extends to an action
of the Brauer monoid $\Br(M)$ determined by the following rules for the
generators $e_i$, where $i\in \M$ and $B\in \AO$.
\begin{eqnarray*}
\delta B &=& B,\\
e_iB &=& \begin{cases}B&\mbox{if  } \a_i\in B,\\
\cl{(B\cup\{\a_i\})}&\mbox{if }\a_i\perp B,\\
r_{\b}r_i B &\mbox{if }\b\in B\setminus\a_i^\perp.
\end{cases}
\end{eqnarray*}
For this action, if $a\in \BrM(M)$ and $B$, $C\in\AO$ satisfy $B\subseteq
C$, then $aB\subseteq aC$.
\end{Prop}

\np When considering words $\und a$ of $F_n$ and  $B\in\AO$,
we will often write
$\und a B$ instead of $\pi(\und a) B$ and
$B\und a $ instead of $\pi(\und a)^{\op} B$. The latter defines right
actions of $F_n$ and $\BrM(M)$ on $\AO$.

\begin{Remark}
\label{rmk:attack}
\rm Let $M=\D_n$ and suppose that $\und a$ is a word in $F_n$. Then, by
Lemma \ref{lm:Worbits}, $\und{a}(\emptyset) \in W B_Y$ for some $Y\in\MY$.
In \cite[Proposition~$4.9$]{CFW}, it is shown that, up to powers of
$\delta$, the elements $\pi(\und a)$ of $\BrM(\D_n)$ are in bijective
correspondence with triples $(B,B',z)$ with ${\und a}(\emptyset)$ and ${\und
a}^\op(\emptyset)$ in the same $W$-orbit $WB_Y$ for $Y\in\MY$ and $z\in
W(M_Y)$, where $M_Y$ is the Coxeter type corresponding to $Y$ as specified
in column 3 of Table \ref{table:types2}.  We will prove a counterpart of
this result.  In fact, we will prove a stronger statement
(Corollary~\ref{cor:finalreduction}) about rewrites of $\und a$ in $F_n$
rather than equality for $\pi(\und a)$ in $\BrM(\D_n)$.  In general, the
expressions appearing in [loc.~cit.]  are not reduced, which makes them
unsuitable for rewrite purposes.  By way of example, we mention that, for
$e_nr_n^*$, we will find an expression of height $1$ (namely $\und z_n^*$ as in
Lemma \ref{lm:enrn*}) where
$r_n^*$ (see Definition~\ref{def:rn*}) has large height.
\end{Remark}

\section{Elementary properties of the Brauer monoid}\label{sec:reduction}
In this section we prepare for the reduction of words in $F_n$ for the
Brauer monoid, $\BrM(\D_n)$, of type $\D_n$. The purpose of this and the
next section is to show that, up to homogeneous equivalence, each element of
$F_n$ has a unique reduced word. This goal is achieved in Corollary
\ref{cor:finalreduction}.

\np We will use the action of Proposition \ref{prop:BrMonAction}.
Let $\und{a}\in F_n$.

Our immediate goal will be to show that $\und{a}$ can be rewritten to a
reduced word that is uniquely determined by $\pi(\und a)$ up to homogeneous
equivalence, so the reduced word will be a unique element of $F_n/\isog$.  In
fact, we shall be working with words in $F_n$ but often think of them as
representing classes in $F_n/\isog$.

Later, in Sections~\ref{sec:conclusion} and \ref{cellular}, we will use
words in $F_n$ to represent monomials in $\BMW(\D_n)$.
Before we continue we introduce some notation.

\begin{Notation}
\label{not:support}
\rm
Suppose that $k$ and $i$ are
two nodes of $\D_n$.  Let $i=i_1,i_2,\ldots,i_r=k$ be the geodesic path from
$i$ to $k$ in $\D_n$. Then we set $e_{i,k} = e_{i_1}e_{i_2}\cdots
e_{i_r}$, which we interpret as an element of $F_n$.  Notice the first factor
is $e_i$ and the last is $e_k$.  In particular, for $i<k$ and $k\ge3$, we have
$e_{i,k}=e_ie_{i+1}\cdots e_k$ unless $i=1$ in which case it is
$e_1e_3e_4\cdots e_k$.  Also $e_{1,2}=e_1e_3e_2$ is a special case.

Let $\b$ be a positive root.  If $\b=\sum \lambda_i \a_i$ we call the
\emph{support} of $\b$ the set of nodes $i$ for which $\lambda_i\ne 0$; it is
denoted $\Supp(\b)$. As in \cite{CGW}, we will write, if $k$ is a node of the
diagram, $\Proj(k,\b)$ for the node of $\D_n$ in $\Supp(\b)$ nearest to $k$.
There is a unique one as the support is a connected set of nodes in the
Dynkin diagram $\D_n$, which is a tree.
\end{Notation}

\begin{Def}\label{df:abetan}
\rm If $k\in\Supp(\b)$, then, as follows directly from \cite[Proposition
  3.2]{CGW}, there is a unique Weyl group element $a_{\b,k}$ of smallest
length that maps $\{\a_k\}$ to $\{\b\}$ in the action of Proposition
\ref{prop:BrMonAction} (so $a_{\b,k}\{\a_k\}=\{\b\}$).  Its height, as a
monomial of $\Br(\D_n)$, is equal to $\het(\b)-1$.  The opposite element
$a_{\b,k}^{\op}$ maps $\{\b\}$ to $\{\a_k \}$.  We will often view
$a_{\b,k}$ as an element of $F_n$ in the guise of a shortest expression for
$a_{\b,k}$ as a product of simple reflections. Since any two such
expressions are homogeneously equivalent, they represent the same element of
$F_n/\isog$, which suffices for our purpose of reductions.

We extend the definition of
$a_{\b,k}$ to the case where $k\not\in\Supp(\b)$.
For $\b$ a positive root with $k\not\in\Supp(\b)$ and $k'$ the node
next to $k$ on the geodesic path from $k$ to $j=\Proj(k,\b)$,
we set $a_{\b,k}=a_{\b,j}e_{j,k'}$ in $F_n$.
\end{Def}

\np
We will be mainly concerned with the case $k=n$.

\begin{Lm}\label{abetan}
The elements $a_{\b,n}$ satisfy the following properties.
\begin{enumerate} [(i)]
\item If $j\leq n-1$, then  $a_{\a_j,n}e_n = e_{j,n}$.
\item If $j$ is a node of $\D_n$ such that $\b-\a_j$ is a root,
then $a_{\b,n}e_n \isog r_ja_{\b-\a_j,n}e_n$.
\item $\het(a_{\b,n}) = \het(\b) -1 $.
\end{enumerate}
\end{Lm}

\nl\begin{proof}
(i). Clearly $n$ is not in the support of $\a_j$ and so
$a_{\a_j,n}=e_{j,n-1}$.
The required equality follows from multiplication by $e_n$ on the right.

\nl(ii).  We first consider the case where $n\in\Supp(\b)$.  In this case
$a_{\b,n}$ is any word of shortest length which takes $\a_n$ to $\b$.  Its
length is $\het(\b)-1$.  As mentioned above and in \cite[Proposition 2.3]{CGW}
it is a unique up to homogeneous equivalence.  If $j\neq n$, then
$a_{\b-\a_j,\a_n}$ is a word of shortest length taking $\a_n$ to $\b-\a_j$ and
so $r_ja_{\b-\a_j,\a_n}$ is a word of shortest length taking $\a_n$ to $\b$,
proving that $ a_{\b,n}\isog r_ja_{\b-\a_j,n}$, and
so $a_{\b,n}e_n\isog r_ja_{\b-\a_j,n}e_n$.

If $j=n$, then $\b-\a_n$ is a root.  Because of the structure of the roots
of $\D_n$, this means the coefficient in $\b$ of both $\a_n$ and $\a_{n-1}$
as a linear combination of simple roots is $1$ and so $\b-\a_n$ has ${n-1}$
in its support but not $n$.  In particular, $a_{\b-\a_n,n-1}$ is a word of
height $\het(\b)-1$ taking $\a_{n-1}$ to $\b-\a_n $.  As $n$ is not in the
support of $\b-\a_n$ but ${n-1} $ is, $a_{\b-\a_n,n-1}$ is a word in $r_i$
with $i\leq n-1$ and further as the coefficient of $\a_{n-1}$ in $\b-\a_n$
is just $1$, all the $r_i$ occurring in a reduced word for $a_{\b-\a_n,n-1}$
have $i\le n-2$.  In particular $r_n$ and $a_{\b-\a_n,n-1}$ commute.  Also,
$a_{\b-\a_n,n}=a_{\b-\a_n,n-1}e_{n-1}$ by definition.  Now
$r_na_{\b-\a_n,n}e_n=r_na_{\b-\a_n,n-1}e_{n-1}e_n\isog
a_{\b-\a_n,n-1}r_ne_{n-1}e_n$.  By (HNree) $r_ne_{n-1}e_n\isog r_{n-1}e_n$.
In terms of the action of Proposition \ref{prop:BrMonAction}, this implies
$a_{\b-\a_n,n-1}r_{n-1}\{\a_n\}=a_{\b-\a_n,n-1}\{\a_{n-1}+\a_n\}$.  Recall
$a_{\b-\a_n,n-1}$ is a shortest word in $r_1,\ldots,r_{n-1}$ taking
$\{\a_{n-1}\}$ to $\{\b-\a_n\}$ and so is a word of shortest length taking
$\{\a_{n-1}+\a_n\}$ to $\{\b\}$ as $a_{\b-\a_n,n-1}$ fixes $\{\a_n\}$.  Now
$a_{\b-\a_n,n-1}r_{n-1}$ is a shortest word taking $\{\a_n\}$ to $\{\b\}$
and so $a_{\b-\a_n,n-1}r_{n-1}\isog a_{\b,n}$.  This gives
$r_na_{\b-\a_n,n}e_n\isog a_{\b-\a_n,n-1}r_{n-1}e_n \isog a_{\b,n}e_n$.

If $n\notin \Supp\ (\b)$, let $i=\Proj(\b,n)$.  If $i>3$, the argument above
applies directly with $i$ instead of $n$ and $j\le i$, giving
$r_ja_{\b-\a_j,i}e_i\isog a_{\b,i}e_i$. The assertion now follows from right
multiplication by $e_{i+1,n}$.  For $i=3$, the root $ \b$ is $\a_3+\a_2$ or
$\a_3+\a_1$ and the arguments are similar.  Notice $i$ cannot be $1$ or $2$
as $\b-\a_j$ is a root.

\nl(iii). This is direct from (ii) and the definition of $a_{\b,n}$.
\end{proof}

\np
\begin{Remark}\label{rmk:notInBMW} \rm
As the proof uses the relation (HNree) which is not binomial in the BMW
algebra, two homogeneously equivalent words of (ii) do not necessarily have
the same image under $\rho$ in the BMW algebra. Indeed, if $j=n$ and
$\b=\a_{n-1}+\a_n$, then $a_{\b,n} = r_{n-1}$ and $a_{\b-\a_n,n}=e_{n-1}$, so
$\rho(a_{\b,n}e_n) = g_{n-1}e_n$ is distinct from
$\rho(r_na_{\b-\a_n,n}e_n)=g_ne_{n-1}e_n$. As indicated in
Proposition \ref{prop:wordssameinbrauer},
the two
expressions are equal up to sums of monomials of lower height (with
coefficients in the ideal generated by $m$).
\end{Remark}

\np We have denoted words in $F_n$ by underlined symbols like $\und a$.
In the remainder of the paper we will need to reduce words which have
specific $r_i$ or $e_i$ in them.  It is notationally awkward to have long
strings underlined, and so we will dispense with this for words including
such $r_i$ and $e_i$. For example we write $\und a r_ir_je_i\homog \und
ae_je_i$ rather than $\und{ar_ir_je_i}\homog \und{ae_je_i}.$ We continue to
underline general elements of $F_n$ as $\und a$.

Let $M$ be a Coxeter diagram with $n$ nodes.  The \emph{Matsumoto--Tits
rewrite rules} of type $M$ on $s_1,\ldots,s_k$ are the following rewrite
rules in the free monoid on $s_1,\ldots,s_n$.
\begin{eqnarray*}
s_is_i&\homog& 1\\
s_is_j&\homog& s_js_i\mbox{\ \ \rm if \ \ } i\not\sim j\\
s_is_js_i&\homog& s_js_is_j\mbox{\ \ \rm if \ \ } i\sim j
\end{eqnarray*}
Note that the second and the third rule are homogeneous.

\begin{Lm}\label{lm:Matsumoto}
Let $M$ be a Coxeter diagram with $n$ nodes.  Then any two reduced words
with respect to the Matsumoto--Tits rewrite rules of type $M$ on
$s_1,\ldots, s_n$ are homogeneously equivalent, that is, can be rewritten
into each other by means of a series of the second and the third rewrite
rules.
\end{Lm}

\nl\begin{proof} The result can be found in \cite{Tits69} and is independently
proved in \cite{Matsum}.  A more general version is found in \cite{Bour}.

\end{proof}

\np As a first application, note that, for the subgroup $W$ of $\BrM(\D_n)$,
the rewrite rules with $r_1,\ldots,r_n$ instead of $s_1,\ldots,s_k$ coincide
with (RSrr), (HCrr), and (HNrrr) of Table \ref{BrauerTable}. Therefore, each
element of $W$ corresponds to a unique reduced word of $F_n$ up to homogeneous
equivalence.  In other words, the equivalence classes in $F_n$ of reduced
words over $\{r_1,\ldots,r_n\}$ correspond bijectively with the elements of
the Coxeter group $W$.  This implies that, for each reduced word $\und a\in
F_n$ all of whose symbols are in $\{r_1,\ldots,r_n\}$, its homogeneous
equivalence class is uniquely determined by $\pi(\und a)$.
In Proposition \ref{prop:ZY}, we will
generalize this application to recognize Coxeter groups of type $M_Y$
for each $Y\in\MY$,
using words $\und{s}_i$ in $F_n$ to be specified in Notation
\ref{not:zn}.

A slightly less general statement holds for $\BMW(\D_n)$ instead of
$\BrM(\D_n)$. As of the above-mentioned rewrite rules, (HCrr) and (HNrrr)
are binomial in Table \ref{BMWTable} as well, for each reduced $\und a\in
F_n$ all of whose symbols are in $\{r_1,\ldots,r_n\}$, its homogeneous
equivalence class is uniquely determined by $\rho(\und a)$ as well.  In
Proposition \ref{prop:hecke}, we will generalize this application, using the
same words $\und{s}_i$ as above in $F_n$ to recognize subquotients
of $\BMW(\D_n)$ isomorphic to Hecke algebras of type $M_Y$ for $Y\in \MY$.

Observe that $F_{n-1}$ is a submonoid of $F_{n}$.

\begin{Lm}\label{lm:enrn*} Let
$\und{z}_n^* = e_{n,2}r_1e_{3,n}$ for $n\geq 3$, $\und{z}_1^*=r_2e_1$, and $\und{z}_2^*=r_1e_2$  all of these viewed as words in $F_n$.
Then $\und{z}_n^*$ has height $1$
and occurs in the following reductions for $n\ge3$.
\begin{enumerate}[(i)]
\item
$r_n^*e_n\homog \und{z}_n^*$ and $e_nr_n^*\homog \und{z}_n^*$.
\item
$\und{z}_n^* \isog e_{n,3}r_2e_1e_{3,n}$.
\item
\label{lm:xin}
For $n\ge4$ and
$i\in\{1,\ldots,n-2\}$, $e_i\und{z}_n^* \isog \und{z}_i^* e_n \isog e_n\und{z}_i^*$ and $r_i\und{z}_n^*
\isog \und{z}_n^*r_i$.
\item
$\und{z}_n^*e_{n-2}
\isog e_n\und{z}_{n-2}^*$ and $e_n\und{z}_n^*\isog \und{z}_n^*e_n\isog \delta \und{z}_n^*$.
\item
\label{en*sqared}
$\und{z}_n^*\und{z}_n^*\homog
\delta e_n$.
\end{enumerate}
For $n$ equal to $1$ or $2$, statements (i) and (v) also hold.
\end{Lm}

\nl\begin{proof} Assume first $n\geq 3$.  By definition, there is only one
factor $r_i$ in $\und{z}_n^*$ and so its height is at most $1$.

To see that it does not have height zero we use the representation
$\rho_{W\{a_n\}}$ of \cite[Theorem 3.6]{CFW}. In particular we are
considering $Y(1)=\{ \a_n\}$. Consider the action of $\und z_n^*$ in the
notation of [loc.~cit.] on the 1-space spanned by the vector
$\xi_{\{\a_n\}}$.  Indeed $e_{2,n}\xi_{\{\a_n\}}=\delta\xi_{\{\a_2\}}$ and
then $r_1\xi_{\{\a_2\}}=\xi_{\{\a_2\}}h_{1,\a_2}$.  It follows from
[loc.~cit.] that $h_{1,\a_2}$ is one of the generators of $M_{\{n\}}$ which in
this case is the Weyl group of the diagram of type $\A_1\D_{n-2}$.  Now act
by $e_{3,n}$ to see $\und
z_n^*\xi_{\{\a_n\}}=\delta\xi_{\{\a_n\}}h_{1,\a_2}$.  If $\und z_n^*$ could
be reduced it would have height~$0$ and the action on $\xi_{\{\a_n\}}$ would
either be~$0$ or would be $\xi_{\{\b\}}\delta^k$ for $\b$ a root and for
some $k$ a contradiction.  This means $\und z_n^*$ has height $1$.

\nl(i).
Let
$w_{2,n} =
r_3r_2r_4r_3r_5r_4\cdots r_{n-1}r_{n-2}r_nr_{n-1}$
be as in
\cite[Lemma 3.1]{CGW}
and set $w_{n,2} = w_{2,n}^\op$. Then
$r_n^* = w_{n,2}r_1w_{2,n}$ where $r_n^*$ was defined in Definition~\ref{def:rn*}.
In order to show the required reductions,
we use repeatedly the reducing relation (RNrre), that is, $r_jr_ie_j\homog
e_ie_j$ for $i\sim j$.
In particular, $w_{n,2}e_n\homog e_{2,n}$.  Now
$r_1e_{2,n}\homog e_2r_1e_{3,n}$ and $r_n^*e_n = w_{n,2}r_1w_{2,n}e_n\homog
e_{n,2}r_1e_{3,n}=\und{z}_n^*$.  A similar computation shows
that $e_nr_n^* \homog \und{z}_n^*$.

\nl(ii).  This statement holds because of $e_3e_2r_1e_3 \isog e_3r_2e_1e_3$,
which is immediate from the defining relation (HTeere).

\nl(iii).
For $i\in\{2,\ldots,n-2\}$, by the definition of $e_{k,n}$, (HCee), and
(HNeee),
\begin{eqnarray*}
e_i\und{z}_n^* &\isog &  e_{n,i+2} e_ie_{i+1}e_{i,2}r_1e_{3,n}
  \isog   e_{n,i+2}e_{i,2}r_1e_{3,n} \\
   &\isog &  e_{n,i+2} e_{i,2} r_1e_{3,i}    e_{i+1}   e_{i+2,n}
  \isog   e_{i,2} r_1e_{3,i} e_{n,i+2}e_{i+1}  e_{i+2,n}        \\
    &\isog &    \und{z}_i^*   e_n
\isog    e_n\und{z}_i^*  .
\end{eqnarray*}

By (HCer), (HNree), and (HNeer),
\begin{eqnarray*}
r_i\und{z}_n^* &\isog &  e_{n,i+2} r_ie_{i+1}e_{i,2}r_1e_{3,n}
\isog   e_{n,i+2}r_{i+1}e_ie_{i-1,2} r_1e_{3,n} \\
& \isog &  e_{n,i+2}e_{i+1}r_{i+2}e_{i,2}r_1e_3e_{4,n}
\isog   e_{n,2}r_1e_{3,i}r_{i+2}e_{i+1}e_{i+2,n}  \\
&\isog & e_{n,2}r_1e_{3,i}r_{i+1}e_{i+2,n}
\isog   e_{n,2}r_1e_{3,i}e_{i+1}r_{i}e_{i+2,n}  \\
&\isog &  \und{z}_n^* r_i.
\end{eqnarray*}

The case $i=1$ is notationally different but can be done the same way
as $i=2$.

\nl(iv).  In view of the palindromic nature of the word $\und{z}_i^*$ and the fact,
proved in (iii), that $\und{z}_i^*$ and $ e_n$ commute homogeneously, we see
that $e_i $ and $\und{z}_n^*$ commute homogeneously. Applying this with
$i=n-2$ gives $\und{z}_{n-2}^*e_n\isog e_n\und{z}_{n-2}^*$.
The second chain of homogeneous equivalences is a direct consequence of (RSee).

\nl(v).  By (RSee), (HCer), (HNeee), and (RSrr),
\begin{eqnarray*}
\und{z}_3^*\und{z}_3^* &=&  e_3e_2r_1e_3e_3e_2r_1e_3
 \homog  \delta e_3r_1e_2e_3e_2r_1e_3
                                   \homog  \delta e_3r_1e_2r_1e_3   \\
                                   &\homog & \delta e_3 r_1 r_1 e_2 e_3
                                    \homog  \delta e_3e_2e_3
                                  \homog   \delta e_3.
\end{eqnarray*}
Also, by (HSee), (HCer), (HNeee), and (RSrr),
\begin{eqnarray*}
\und{z}_n^*\und{z}_n^* &= & e_{n,4}e_3r_1e_2e_3e_{4,n}e_{n,4}e_3r_1e_2e_3e_{4,n} \homog  \delta e_{n,4}e_3r_1e_2e_3e_2r_1e_3e_{n,4} \\
                                       & \homog & \delta e_{n,4}e_3e_{4,n}
                                       \homog   \delta e_n.
\end{eqnarray*}
The cases $n=1$ and $n=2$ can be done separately.
\end{proof}

\begin{Notation}
\label{not:zn}
\rm
Let $M$ be of type $\ADE$.
For any coclique $Y$ of $M$, we write
$e_Y = \prod_{y\in Y} e_y$ and
$\he_Y = \delta^{-|Y|}\prod_{y\in Y} e_y$.
All factors commute, so we need
not care about the order in which they occur.
For instance $Y(0)=\emptyset$ and $\he_{Y(0)}=1$, whereas
$Y(1) = \{\a_n\}$ and $\he_{Y(1)}=\he_n$.

\rm
We distinguish the following elements of $F_n$ according to the different
possibilities for $Y\in\MY$. We need
$\und{z}_n^*$ as in Lemma \ref{lm:enrn*} and
$e_n^* = e_{n,2}e_1e_{3,n}$, which is the height zero analog of $\und z_n^*$.
The elements ${\und s}_i$ and ${\und f}_i$ will play roles reminiscent of
$r_i$ and $e_i$.
\begin{eqnarray*}
Y=Y(t) \hfil \qquad (t>0)&: \ & \und{s}_0 =
 \und{z}_{n}^* \delta^{-1}\he_Y,
\  \und{s}_i = r_i\he_Y, \\
&&
 \und{f}_0 =
 \und{e}_{n}^* \delta^{-1}\he_Y,
\  \und{f}_i = e_i\he_Y \ \ (1\le i\le n-2t)\\
Y=Y(0)\hfill&:\ &  \und{s}_i = r_i, \ \und{f}_i = e_i \ \
(1\le i\le n)\\
Y=Y(\frac{n-1}{2})\hfill \qquad(n \mbox{ odd})&: \ & \und{s}_0 =
 \und{z}_{n}^* \delta^{-1}\he_Y,
\ \und{f}_0 = \und{e}_{n}^* \delta^{-1}\he_Y
\\
Y=Y(\frac{n}{2})\hfil \qquad(n \mbox{ even})&: \ & \und{s}_1 = r_{1}^*
 \delta^{-1}\he_Y,\
 \und{f}_1 = \und{\he}_{2} \he_Y
\\
Y=Y'(\frac{n}{2})\hfil \qquad(n \mbox{ even})&: \ & \und{s}_1 = r_2^* \delta^{-1}\he_Y,
\ \und{f}_1 = \und{\he}_{1} \he_Y
\\
Y=Y^*(t)\hfill\qquad (t>0)&:\ &  \und{s}_i = r_{i+3}\he_Y, \\
&&
 \und{f}_i =
 \und{e}_{i+3} \delta^{-1}\he_Y \ \ (1\le i\le n-2t-1)
\end{eqnarray*}
\end{Notation}

\np Let $Y\in\MY$. The indices $i$ of $\und{s}_i$ and $\und{f}_i$ occurring
in Notation \ref{not:zn} are attached to the diagram $M_Y$ in such a way
that $0$ (if it occurs) corresponds to the isolated component $\A_1$ of
$M_Y$ and the other component (of type $\A$ or $\D$) is labeled as usual for
$\A$ and as indicated in Figure \ref{fig:Dn} for $\D$.
For instance, in case $Y=Y(t)$ with $0<t<(n-1)/2$,
the diagram $M_Y = \A_1\D_{n-2t}$ is labeled as follows.
$$\DnS
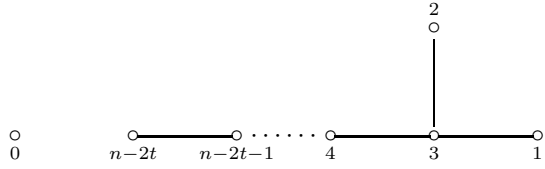
$$

\np In the proposition below we establish the Matsumoto--Tits rewrite rules
for the Coxeter group of type $M_Y$ with generators $\pi({\und s}_i)$ as in
(\ref{not:zn}) and identity $\he_Y$.  The $\und{f}_i$ will be studied in the
next section.

\begin{Prop}\label{prop:ZY}
Let $n\ge4$ and
$Y\in\MY$.
The words $\und{s}_i$ in $F_n$,
for $i$ a node of $M_Y$, have height $1$ and satisfy
the following properties.
\begin{enumerate}[(i)]
\item With respect to the rewrite system of Table \ref{BrauerTable} in
$F_n$, the words $\und{s}_i$ satisfy the
Matsumoto--Tits rewrite rules of type $M_Y$ with identity element
$\he_Y$. That is, they satisfy $\he_Y\he_Y\isog \he_Y$, $\und{s}_i\he_Y\isog
\und{s}_i$, $\he_Y\und{s}_i\isog \und{s}_i$, $\und{s}_i \und{s}_i\homog
\he_Y$, $\und{s}_i\und{s}_j\isog \und{s}_j\und{s}_i$ if $i\not\sim j$, and
$\und{s}_i\und{s}_j\und{s}_i\isog \und{s}_j\und{s}_i\und{s}_j$ if $i\sim j$,
where $i$ and $j$ are nodes of $M_Y$.
\item The elements $\pi(\und{s}_i)$,
for $i$ running through the nodes of $M_Y$,
generate a Coxeter group of type $M_Y$ in $\BrM(\D_n)$ with
identity element $\pi(\he_Y)$.
\item
\label{not:VYt}
For $Y\in\MY$, denote $U_Y$ the set of words in $F_n\he_Y$ that are
minimal expressions in the $\und{s}_i$ (where $i$ runs over the nodes of
$M_Y$) for elements of the Coxeter group of (ii).  Then the restriction of
$\pi$ to $U_Y$ induces a bijection from the set of homogeneous equivalence
classes in $U_Y$ onto this Coxeter group.
\end{enumerate}
\end{Prop}

\nl\begin{proof} Recall $\het(e_Y) = 0$.  By Lemma \ref{lm:enrn*}(ii),
$\het(\und{z}_n^*) = 1$, so $\het(\overline{s}_0) = 1$, and, clearly,
$\het(\und{s}_i) = 1$ for $i>0$.

\nl(i). We verify the individual rewrite rules in the case where $Y=Y(t)$
and leave the other cases to the reader (as they are similar or easier).
Those involving $\he_{Y(t)}$ at the left hand side are straightforward
applications of the rules (HSee), (HCee), (HCre), and (HNeee).

\nl $\und{s}_i\und{s}_i\homog \he_{Y(t)}$.  By Lemma~\ref{lm:enrn*}(iv) we see
$e_n\und{z}_n^*\isog \und{z}_n^*e_n$ and $\und{z}_n^*e_j\isog e_j\und{z}_n^*$
for $j\leq n-2$ and so $\und{z}_n^*e_{Y(t)}\isog e_{Y(t)}\und{z}_n^*$.  Hence
$\und{z}_n^*$ commutes homogeneously with $e_{Y(t)}$ and so in the definition
of $\ov {s}_0$ it does not matter on which side $e_{Y(t)}$ occurs. In
particular, using Lemma \ref{lm:enrn*}(v), we find $\und{s}_0\und{s}_0\isog
\und{z}_n^*\und{z}_n^*e_{Y(t)}e_{Y(t)} \delta^{-2-2t}\homog\delta
e_ne_{Y(t)}\delta^{-2-t} \homog e_{Y(t)}\delta^{-t}$, which is the identity
element of $\pi(U_{Y(t)})$.  This settles the case $i = 0$.  For $i>0$, the
assertion $\und{s}_i\und{s}_i\homog \he_{Y(t)}$ follows directly from the fact that
$e_{Y(t)}$ and $r_i$ commute and (HSrr).

\nl $\und{s}_i\und{s}_j\homog \und{s}_j\und{s}_i\mbox{\ \ \rm if \ \ }
i\not\sim j$.  For $i=0$ and $j>0$, this follows from Lemma
\ref{lm:enrn*}(iii).  For $i>0$ and $j>0$, it is immediate from (HCrr).

\nl $\und{s}_i\und{s}_j\und{s}_i\homog \und{s}_j\und{s}_i\und{s}_j\mbox{\ \
\rm if \ \ } i\sim j$.  Here we must have $i,j>0$. Now it is immediate from
(HNrrr).

\nl(ii).
The fact that the $\pi(\und{s}_i)$ $(0\le i \le n-2t)$
generate a quotient of the Coxeter group of
type $M_{Y(t)}$ is immediate from (i) and the fact that a rewrite rule
$x\homog y$ in $F_n$ implies $\pi(x) = \pi(y)$.
Therefore, it suffices to show that there is a surjective homomorphism
from the group generated by the  $\pi(\und{s}_i)$ onto
$W(M_{Y(t)})$.
This follows from \cite[Lemma 1.3]{CFW}.

\nl(iii).
This is immediate from (ii) and Lemma \ref{lm:Matsumoto}.
\end{proof}

\section{Reduction in the Brauer monoid}
\label{sec:reduction2}
In this section we continue to discuss reductions of words in $F_n$.  The
main purpose is to show that each word in $F_n$ can be reduced to a
particular form described in Theorem \ref{thm:notthetaunique}.  We first
study the product of a generator and a word $a_{\b,n}e_n$, which is in
reduced form.  Here $a_{\b,n} $ is given by Definition~\ref{df:abetan}.  In
the action of Proposition \ref{prop:BrMonAction}, the element
$\pi(a_{\b,n}e_n)$ maps $\emptyset$ to $\{\b\}$, so after left
multiplication by $e_i$ it will map $\emptyset$ to $\{\a_i\}$ or (in case
$\a_i\perp \b$) to $\{\a_i,\b\}$, and, after left multiplication with $r_i$,
it will map $\emptyset$ to $\{r_i\b\}$. The lemma below will find
corresponding reduced words.  In order to control the kernel of this action,
we need a little more notation.

\begin{Notation}\label{not:ZY}
\rm For $Y\in\MY$, let $Z_Y$ be the subsemigroup of $F_n$ generated by all
$\delta^{j}e_Y$ for all $j\in\Z$, and $\und{s}_i$ and $\und{f}_i$ for all
nodes $i$ of $M_Y$ as in Notation \ref{not:zn}.
We also write $Z_n$ instead of $Z_{Y(1)}$.
The subsemigroup
$Z_{\emptyset}$ coincides with $F_n$.
\end{Notation}

\begin{Lm}\label{lm:oneRootDn}
Let $M = \D_n$, let $i\in\{1,\ldots,n\}$, and let $\b\in\Phi^+$.
Then the word
$e_ia_{\b,n}e_n$
reduces to a word in
$a_{\b',n}Z_n$, where $\b'$ is a positive root with $\het(\b')\le\het
(\b)$.
Also, $r_ia_{\b,n}e_n$
can be reduced to a
word in $a_{\b',n}Z_n$, where $\{\b'\} = r_i\{\b\}$.
Moreover, if $\und a\in F_n$,
then $\und a e_n$ can be reduced to a word in $a_{\b',n}Z_n$, where $\b'$ is
a positive root with $\het(\b')\le \het(\und a)$ and $\b'\in
\und{a}\{\a_n\}$.
\end{Lm}

\nl\begin{proof}
We proceed by induction on $\het(\b)$.
If $\het(\b)=1$ we have $\b=\a_j$ for some
node $j$ of $\D_n$.
By Lemma~\ref{abetan}, $a_{\b,n}e_n\homog e_{j,n}$.

Consider first $e_ia_{\b,n}e_n$.  By the above, $e_ia_{\b,n}e_n\homog
e_ie_{j,n}$.  If $\{i,j\}= \{1,2\}$, then $ e_i e_{j,n} = e_1 e_{2,n} \homog
\delta^{-1} e_2 e_{3,n}e_n^* = \delta^{-1} a_{\a_2,n}e_ne_n^* =
a_{\a_2,n}{\und f}_0 \in a_{\a_2,n}Z_n$.  By symmetry of the diagram, the
case $j=1$ can be replaced by $j=2$ and handled in a similar way, so assume
$j\ge 2$.  If $i<j$, then $e_i$ can be commuted to the right and be absorbed
into $Z_n$ as $\und f_i$.  If $i=j-1$, then we may assume $i\ge2$ as we
already handled the case $\{i,j\} = \{1,2\}$, and so $e_ia_{\b,n}e_n \homog
e_ie_{j,n}=e_{i,n}=a_{\a_i,n}e_n$.  If $i=j$ we obtain $e_ia_{\b,n}e_n \homog \delta
e_{i,n}$.  If $i=j+1$ we can use $e_ie_je_i\isog e_i$ to derive
$e_ia_{\b,n}e_n \homog e_{i,n}$.  Otherwise
$i\not\sim j$ and $i\not= j$; commute the $e_i$ past terms in
$e_{j,n}$ to obtain $e_ia_{\b,n}e_n \homog e_je_{j+1}\cdots
e_ie_{i-1}e_i\cdots e_n$.  Now use $e_ie_{i-1}e_i\isog e_i$ and commute the
preceding terms $e_k$ to the right and absorb them into $Z_n$ as products of $\und f_k$.
In each of these cases $e_ia_{\b,n}e_n \homog a_{\b',n}\und z$ for some
$\und z\in
Z_n$ and some $\b'\in \{\a_i,\a_j\}=e_i\{\a_j\}$, as required.

We now consider $r_ia_{\b,n}e_n$ with $\b=\a_j$ for some node $j$, where
$r_ia_{\b,n}e_n\homog r_ie_{j,n}$.  There are two special cases which we
handle directly, viz., $j=2$ with $i=1$ and $j=1$ with $i=2$.  For the first
we have $r_1e_{2,n}\homog e_2r_1e_{2,n} \delta^{-1}$.  Notice $e_2\homog
e_{2,n}e_{n,2}\delta^{-1}$ and so $r_1e_{2,n}\homog e_2r_1e_{3,n}\homog
e_{2,n}e_{n,2}r_1e_{3,n} \delta^{-1} = e_{2,n}\und{z}_n^* \delta^{-1} =
a_{\a_2,n}\und{z}_n^*\delta^{-1}$ and we are done as
$\und{z}_n^*\delta^{-1}=\und s_0\in Z_n$.  The other case is similar.
Assume, therefore, that these special cases do not occur.  If $i=j-1$, we
have $r_ie_{j,n-1}e_n = a_{r_i\b}e_n$ and we are done.  If $i<j-1$, then
$r_i$ commutes homogeneously through to give $e_{j,n}r_i$, unless we have
$i=1$ and $j=3$, a case that can be treated as $i=2$ and $j=3$, which is
done below; observe that the expression $e_{j,n}r_i$ is equal to $
e_{j,n}\und {s_i}=a_{a_j,n}e_n\und{s_i}$ and satisfies all the requirements.
If $i=j$, use $r_ie_i\homog e_i$ to see that $r_ie_{j,n-1}e_n \homog
e_{i,n}\in a_{\b,n}Z_n$.  As above if $i=j+1$, then by (HNree),
$r_{j+1}e_je_{j+1}e_{j+2,n}\homog
r_je_{j+1}e_{j+2,n}=r_je_{j+1,n}=a_{\a_j+\a_{j+1},n}e_n$.  This is what is
required as here $\b=\a_j$ and $r_{j+1}\a_j=\a_j+\a_{j+1}$.  Otherwise,
$i>j+1$ and $r_ie_{j,n}\homog e_{j,i-2}r_ie_{i-1}e_{i,n}\homog
e_{j,i-2}r_{i-1}e_{i,n}$.  Now if $j<i-2$ we see
$e_{j,i-2}=e_{j,i-3}e_{i-2}$, and we use $e_{i-2}r_{i-1}\homog
e_{i-2}e_{i-1}r_{i-2}$ to find $e_{j,i-1}r_{i-2}e_{i,n}$ and, commuting
$r_{i-2}$ homogeneously to the right, we obtain the required form.

We may suppose then that $\b$ has height greater than~$1$ and so there is a
node $j$ for which $\b-\a_j$ is a root.  Throughout this part of the proof
we use Lemma~\ref{abetan} when $\b-\a_j$ is a root to see that
up to homogeneous equivalence $a_{\b,n}e_n$ $\isog$ $r_ja_{\b-\a_j}e_n$.

Again, consider first $e_ia_{\b,n}e_n$.  Choose $j=i$ if possible.  If so,
we use $e_ir_i\homog e_i$ to obtain $e_ir_ia_{\b-\a_i,n}e_n\homog
e_ia_{\b-\a_i,n}e_n$.  The resulting word has lower height than
$e_ir_ia_{\b-\a_i,n}e_n$ and we use induction to finish.  Suppose $i \not
\sim j$ and $i\ne j$.  Then $e_ir_ja_{\b-\a_j,n}e_n\homog
r_je_ia_{\b-\a_j,n}e_n$.  Now apply the induction hypothesis to
$e_ia_{\b-\a_j,n}e_n$ so $e_ia_{\b-\a_j,n}e_n\homog a_{\b',n}e_n\und z$
where $\und z\in Z_n$ and $\het(\b')< \het (\b)$.  In view of this
inequality, induction applies to the statement involving $r_ja_{\b',n}e_n$.
Acting by $r_j$ could raise the height at most one, still leaving $\het
(r_j\b')\leq \het(\b)$ as needed.  Suppose $i\sim j$.  We know that
$(\a_i,\b) $ is not $1$ as we have chosen $j=i$ if possible above.  This
means either $(\b,\a_i)=0$ or $(\b,\a_i)=-1$.  Suppose first $(\b,\a_i)=0$.
Then $(\b-\a_j,\a_i)=1$ and so $\b-\a_j-\a_i$ is a root and
$a_{\b,n}e_n=r_jr_ia_{\b-a_j-\a_i,n}e_n$; now
$e_ir_jr_ia_{\b-\a_j-\a_i,n}e_n\homog e_ie_ja_{\b-\a_j-\a_i,n}e_n$, and we
can finish by induction to get the result as the height of the root
$\b-\a_j-\a_i$ is at most $\het (\b)-2$.  Suppose now $(\b,\a_i)=-1$.  Then
$e_ia_{\b,n}e_n \homog e_ir_ja_{\b-\a_j,n}e_n \homog
e_ie_jr_ia_{\b-\a_j,n}e_n $.  Notice $(\b-\a_j,\a_i)=-1+1=0$ and so
$r_i(\b-\a_j)=\b-\a_j$, from which we derive $e_ia_{\b,n}e_n \homog
e_ie_jr_ia_{\b-\a_j,n}e_n \homog e_ie_ja_{\b-\a_j,n}e_n\und z$
for some $\und z \in Z_n$ by the
induction hypothesis for the action of $r_i$.  Using the induction
hypothesis twice more, we find $e_ia_{\b,n}e_n \homog e_ia_{\b',n}e_n \und
z'\homog a_{\b''}e_n\und z''$ for certain roots $\b'$ and $\b''$ whose
height is at most $\het (\b)-1$ and $\und z'$, $\und z''\in Z_n$.  This ends
the part of the proof involving left multiplication by $e_i$.

We now consider $r_ia_{\b,n}e_n$ where $\het (\b)>1$.
If $(\b,\a_i)=-1$, then $r_ia_{\b,n}e_n\isog
a_{\b+\a_i,n}e_n $ by Lemma~\ref{abetan} and we are done.  Suppose
$(\b,\a_i)=1$.  Then $\b-\a_i$ is a root and $a_{\b,n}e_n \isog
r_ia_{\b-\a_i,n}e_n$.  Now use $r_ir_ia_{\b-\a_i,n}e_n\homog a_{\b-\a_i,n}e_n
\isog a_{r_i\b,n}e_n$ to finish.

Therefore, we can assume $(\b,\a_i)=0$.  There is a node $j$ for which
$\b-\a_j$ is a root and so $r_ia_{\b,n}e_n\isog r_ir_ja_{\b-\a_j,n}e_n$ by
Lemma~\ref{abetan}.  The arguments here are similar to the ones at the
beginning of this proof when $\het \b >1$.  In particular, if $i\not \sim j$
and $i\ne j$ this reduces to $r_jr_ia_{\b-\a_j,n}e_n$ and we use induction
for $r_i$ acting in the case $(\a_i,\b-\a_j)=0$.

The only remaining case is $i\sim j$ and still $(\b,\a_i)=0$.  Here
$\b-\a_i-\a_j$ is a root orthogonal to $\a_j$ and $a_{\b,n}e_n\isog
r_jr_ia_{\b-\a_j-\a_i,n}e_n$ by Lemma~\ref{abetan}.  We consider
$r_ir_jr_ia_{\b-\a_j-\a_i,n}e_n$ and so use the homogeneous relation
$r_ir_jr_i\isog r_jr_ir_j$, the induction hypothesis and
$(\a_j,\b-\a_j-\a_i)=0$ to derive $ r_ir_jr_ia_{\b-\a_j-\a_i,n}e_n\homog
r_jr_i a_{\b-\a_j-\a_i,n}e_n\und z \homog a_{\b,n}e_n\und z'$ with $\und z$,
$\und z'\in Z_n$, as required.  This proves all but the last part of the
lemma.

As for the last statement, without loss of generality, we may assume that
$\und ae_n$ is reduced.  We argue by induction on the length of $\und
a$. Whenever $\und a$ is equal to $a_{\b,n}$, there is nothing to show.  In
particular, we may assume that $\und a$ has positive length; say it starts
with $e_i$ or $r_i$. By induction, we have $\und a e_n\homog e_i
a_{\b,n}e_n\und z$ or $\und a e_n\homog r_i a_{\b,n}e_n\und z$ with
$\het(\b)\le \het(\und a)$ for some $\und z\in Z_n$.  The proof now follows
from the second statement in view of $\b' \in a_{\b',n}\{\a_n\}$, which is
clear from the definition of $a_{\b',n}$.
\end{proof}

\np We now return to the sets $Y\in\MY$ and use $Z_Y$ of
Notation~\ref{not:ZY} to reduce words of the form ${\und a}e_Y$.  Sometimes
we come across $e_ne_n^*$, which is homogeneously equivalent to
$e_{n,2}e_1e_2e_{2,n}$ up to powers of $\delta$.  In that case, we usually
invoke Proposition \ref{lm:Theta1} below to reduce the word further. In the
other cases, we have $Y = Y(t)$ for some $t$ or $Y = Y'(n/2)$.

\begin{Notation}
\label{df:Theta}
\rm
By $\Theta$ we denote the ideal
of $\Br(\D_n)$ generated by $e_1e_2$.
For $U$ any subring of
$\Br(\D_n)$, we also write
$Ue_1e_2U$ for the set of all linear combinations of expressions of the form
$ue_1e_2v$ with $u,v\in U$. So
$\Theta = \Br(\D_n)e_1e_2\Br(\D_n)$.
\end{Notation}

\np Note that $Z_{Y^*(t)}$ is contained in $\Theta$ for each $t\in[1,\lfloor
n/2\rfloor]$.

\begin{Prop}\label{lm:Theta1}
Let $Q$ be the subalgebra of $\Br(\D_n)$ generated by all $r_i$ and $e_i$
for $i>1$. Then $Q$ is isomorphic to $\Br(\A_{n-1})$ and satisfies the
following properties.
\begin{enumerate}[(i)]
\item
In $F_n$ any word containing $e_1e_2$ can be reduced to a word of the form
${\und u}e_1e_2{\und v}$ where $\und u$ and $\und v$ are words in
$r_i$ and $e_i$ for
$i>1$, so $\pi({\und u}),\pi({\und v})\in Q$.
\item
The ideal $\Theta$ coincides with $Qe_1e_2Q$. It is isomorphic to the ideal
in $Q$ generated by any $e_i$ $(i>1)$. An explicit height preserving
isomorphism is determined by $u\he_1e_2v\mapsto ue_2v$ for $u,v\in Q$.
\end{enumerate}
\end{Prop}

\nl
\begin{proof} The isomorphism of $Q$ with $\Br(\A_{n-1})$ follows from
the determination of $\Br(\D_n)$ in \cite{CFW}.

\nl(i). This can be shown along the lines of the last paragraph of
\cite[Section 7.1]{CGW}.

\nl(ii).  Let $u,v,u',v'\in Q$.  By considerations in the Brauer algebra of
type $\A_{n-1}$ there is a monomial $h$ in the submonoid of $\BrM(\D_n)$
generated by $e_i$, $r_i$ for $i\ge4$, such that $e_2vu'e_2 = e_2h$.  We
then have $ u\he_1e_2v u'\he_1e_2v' = u\he_1e_2v u'e_2\he_1v' =
u\he_1e_2h\he_1v' = ue_2h\he_1^2v' =  u\he_1e_2hv'$, and the same
multiplication worked out for $ue_2v u'e_2v'$ shows it is equal to
$ue_2hv'$, which proves that the indicated map preserves products. The rank
of domain and range is
$$\sum_{t=1}^{\lfloor n/2\rfloor} \left(
\frac{n!}{2^tt!(n-2t)!}
\right)^2 (n-2t)!
$$
by \cite[Lemma 1.3]{CFW}, and so the map is an isomorphism.
\end{proof}

\np As a consequence, the reduction rules for words mapping into $\Theta$
all follow from reductions in $\BrM(\A_{n-1})$ (applied to elements of the
ideal generated by one of the $e_j$). We will be using these observations
several times below.

The word $e_{Y(t)}$ commutes homogeneously with the
elements $\und{z}_n^*$ and $r_i$, $e_i$ $(i = 1,\ldots, n-2t)$, so, up to
homogeneous equivalence, it does not matter on which side $e_{Y(t)}$ is
located in these expressions for elements of $Z_{Y(t)}$.

\begin{Lm}
\label{lm:mainReduction}
Fix $t\in\{1,\ldots,\lfloor n/2\rfloor\}$.
Consider a word $\und a$ for which $\und a B_{Y(t)}$
is in the same $W$-orbit as $B_{Y(t)}$ and a word
$\und b$ for which $\und b B_{Y^*(t)}$
is in the same $W$-orbit as $B_{Y^*(t)}$.
Then
$\und a e_{Y(t)}$ and  $\und b e_{Y^*(t)}$
each reduce to a word of the form
\begin{eqnarray}\label{eq:abetaks}
&&a_{\b_n,n}a_{\b_{n-2},n-2}\cdots a_{\b_{n-2t+2},n-2t+2}\und z
\end{eqnarray}
with $\b_{n-2k}\in\Phi^+$ for $0\leq k\leq t-1$ such that $\b_{n-2k}$ has
support in $\D_{n-2k}$ for each $ k$, and $\und z\in Z_{Y(t)}$ in the first
case and $\und z\in Z_{Y^*(t)}$ in the second case.  Also, $\b_n\in
\und{a}e_n(\emptyset)$ and $\{\b_n,\b_n^*\} \in \und{b}e_n(\emptyset)$.  The
same applies to $Y'(n/2)$ instead of $Y(n/2)$ if $n$ is even.
\end{Lm}

\nl\begin{proof}
Consider first the case of $\und a$.
The statement that $\b_n\in \und a e_n(\emptyset)$ is straightforward from the
definition and the fact that the terms distinct from $a_{\b_n,n}$
do not move $\a_n=e_n(\emptyset)$.

Notice that
$a_{\b_n,n}a_{\b_{n-2},n-2}\cdots a_{\b_{n-2t+2},n-2t+2}e_{Y(t)}$
is homogeneously equivalent to
$a_{\b_n,n}e_na_{\b_{n-2},n-2}e_{n-2}\cdots a_{\b_{n-2t+2},n-2t+2}e_{n-2t+2}$.
By Lemma \ref{lm:oneRootDn}, $\und a e_n$ can be reduced to
$a_{\b_n,n}e_n \und{z}_n$ for some $\b_n\in\Phi^+$ and $\und{z}_n\in Z_n$.
In particular, up to homogeneous equivalence, cf.\ Lemma
\ref{lm:enrn*}(iii), we may assume $\und{z}_n = \und a'$ or
$\und{z}_n =  \und{z}_n^*    \und {a'}$ for some $\und a' \in F_{n-2}$.
We denote this as $(\und{z}_n^*)^\eps
a'$ where we set $\eps=0$ if it is $\und a'$ and
$\eps=1$ if it is  $\und{z}_n^*    \und {a'}$; so $\eps\in \{0,1\}$.

If $t=1$, we are done by Lemma~\ref{lm:oneRootDn}.  Therefore, we may assume
$t> 1$.  By induction on $n$, we find
$$\und {a'}e_{Y(t)\setminus \{n\}} \homog a_{\b_{n-2},n-2}e_{n-2}\cdots
a_{\b_{n-2t+2},n-2t+2}e_{n-2t+2}\und{z}_{n-2t+2}$$ for some
$\und{z}_{n-2t+2}\in Z_{Y(t)\setminus\{e_n\}}$.  As $\und a\in F_{n-2}$, the
terms of $\pi(\und {a'}e_{Y(t)\setminus \{n\}})$ never include $r_n$,
$r_{n-1}$, $e_n$, or $e_{n-1}$, and so the support of $\beta_{n-2}$ is in
$\D_{n-2}$.

Now, by Lemma \ref{lm:enrn*}, for $\eps\in\{0,1\}$, thanks to $\und
z_n^*e_{n-2}\isog e_n\und z_{n-2}^*\he_{n-2}$, we have, up to powers of
$\delta$
\begin{eqnarray*}
\und a e_{Y(t)}&\homog&
a_{\b_n,n}e_n  (\und{z}_n^*)^\eps    \und {a'} e_{Y(t)\setminus \{n\}}\\
&\homog&
a_{\b_n,n}e_n  (\und{z}_n^*)^\eps
a_{\b_{n-2},n-2}e_{n-2}\cdots a_{\b_{n-2t+2},n-2t+2}e_{n-2t+2}\und{z}_{n-2t+2}\\
&\homog&
a_{\b_n,n}e_n
a_{\b_{n-2},n-2}  e_{n-2}(\und{z}_{n-2}^*)^\eps\cdots a_{\b_{n-2t+2},n-2t+2}e_{n-2t+2}\und{z}_{n-2t+2}\\
&\homog&
a_{\b_n,n}e_n
a_{\b_{n-2},n-2}  e_{n-2}\cdots a_{\b_{n-2t+2},n-2t+2}e_{n-2t+2}
(\und{z}_{n-2t+2}^*)^\eps\und{z}_{n-2t+2}\\
&\homog&
a_{\b_n,n}
a_{\b_{n-2},n-2}  \cdots a_{\b_{n-2t+2},n-2t+2}\und z
\end{eqnarray*}
with $\und z =e_ne_{n-2}\cdots e_{n-2t+2}(\und{z}_{n-2t+2}^*)^\eps\und{z}_{n-2t+2}\in Z_{Y(t)}$
and $\b_{n-2k}$ has support in $\{\a_1,\ldots,a_{n-2k}\}$ for each $k$, as required.

Notice that $\und{z}_{n-2t+2}^* $ homogeneously commutes with elements of
$Z_{Y(t)}$ by Lemma~\ref{lm:enrn*}, (RSer) and (RSre), and the fact that
$\und{z}_{n-2t+2}^*$ starts and ends with $e_{n-2t+2}$.

The case $\und b e_{Y^*(t)}$ runs along the same lines and is simpler in
view of Proposition \ref{lm:Theta1}.
\end{proof}

\np
The special case $t=1$ gives the following corollary.

\begin{Cor}\label{betaimage}
If $\und{a}\{\a_n\}=\{\b\}$ or $\und{b}\{\a_n, \a_n^*\}=\{\b, \b^*\}$, then
$\und{a}e_n\homog a_{\b,n}e_n\und z$ with $\und z\in Z_n$ or
$\und{b}e_n\homog a_{\b,n}e_n\und z$ with $\und z\in Z_{\{\a_n, \a_n^*\}}$.
\end{Cor}

\begin{Remark} \rm
If $B$ is in the $W$-orbit of $B_{Y(t)}$ with $\b\in B$, then $Ba_{\b,n}e_n$
is in the same orbit.  This is clear for the terms from $W$ in $a_{\b,n}$
and also for the terms $e_j$ which map $\a_{j-1}$ to $\a_j$, which is
$r_jr_{j-1}$.  The same is true for $r_3r_1$ moving $\a_1$ to $\a_3$.  Now
$Ba_{\b,n}e_n$ contains $\a_n$ plus roots all in the subsystem of type
$\D_{n-2}$.  The term $\b_{n-2}$ is one of these, which can be associated to
one of the roots of $B$ other than $\b$.  In this way, an order of the roots
of $B$ gives the terms $\b_{n-2s}$ which occur.  The same is true for
$Y'({n}/{2})$ instead of $Y(t)$ if $n$ is even.  A similar result is true
for the case of $Y^*(t)$; here $Ba_{\b,n}e_n$ contains $\a_n$, $\a_n^*$ as
well as roots in $\D_{n-2}$ together with their orthogonal mates.
\end{Remark}

\np We will consider the different ways to write $\und a e_{Y}$ in this reduced
form.  The case of $t=2$ will suffice to argue the
general case.  If $n\ge5$, there are two possibilities, $Y(2)$ and $Y^*(2)$.
As mentioned before, if $n=4$, there
is one more, for $Y = Y'(2)$.

For $Y = Y(2)=\{n, n-2 \}$,
we
consider words of the form $a_{\b_n,n}a_{\b_{n-2},n-2}\und z$ where $\und z\in
Z_{Y(2)}$.  We need a lemma that involves words in $F_n$ mapping $\{\a_n,
\a_{n-2}\}$ to $\{\b, \c\}$ in $\AO$ and the ways to reduce them.
Similarly for $Y = Y^*(2)$ we consider words mapping $\{\a_n, \a_n^*, \a_{n-2}, \a_{n-2}^*\}$ to
$\{ \b,\b^*,\gamma, \gamma^*\}$.

\begin{Lm}\label{lm:twoRoots}
Suppose that $\und a\in F_n$ satisfies $\und a \{\a_n, \a_{n-2}\} = \{\b,
\c\}$ and $\und b\in F_n$ satisfies $\und b \{\a_n, \a_n^*, \a_{n-2},
\a_{n-2}^*\} = \{\b,\b^*, \c,\c^*\}$ .  Then of the two possible reductions
of $\und a e_{Y(2)}$ and $\und b e_{Y^*(2)}$ as in
Lemma~\ref{lm:mainReduction}, at least one can be reduced to the other, that
is, for some $\und z\in Z_{Y(2)}$ or $\und z\in Z_{Y^*(2)}$, respectively,
we have
\begin{eqnarray*}
\mbox{ either }&& \ a_{\b,n}a_{\b_{n-2},n-2}e_ne_{n-2} \homog
a_{\c,n}a_{\c_{n-2},n-2}e_ne_{n-2}\und z\\
\mbox{ or }&& \
a_{\c,n}a_{\c_{n-2},n-2}e_ne_{n-2}
 \homog
a_{\b,n}a_{\b_{n-2},n-2}e_ne_{n-2}\und z.
\end{eqnarray*}
For $n=4$ and $\{\b,\c\}\in WB_{Y'(2)}$,
the same statement holds with $Y'(2)$ instead of $Y(2)$ and
$\{\a_4,\a_1\}$ instead of $\{\a_4,\a_2\}$.
\end{Lm}

\nl
\begin{proof}
We deal with $\und a$ first.
Suppose first that either $\b$ or $\c$ has $n$ in its support.  Without loss
of generality, we assume $n\in\Supp(\b)$.
Then $\pi(a_{\b,n})\in W$ as $n$ is in the support of $\b$.  By (HNeee),
\begin{eqnarray*}
a_{\b,n}e_na_{\b_{n-2},n-2}e_{n-2} &\isog &
a_{\b,n}a_{\b_{n-2},n-2}e_{n-2}e_{n}\\
&\isog & a_{\b,n}a_{\b_{n-2},n-2}e_{n-2}e_{n-1}e_{n-2}e_{n}\\
&\isog &a_{\b,n}a_{\b_{n-2},n}e_{n-2}e_{n-1}e_{n}e_{n-2} \\
    &\isog &a_{\b,n}a_{\b_{n-2},n}e_{n}e_{n-2}\end{eqnarray*}

As $\pi(\a_{\b,n})$ is in the Weyl group,  $a_{\b,n}\{\b_{n-2}\}$ is a single root.  As
$a_{\b,n}a_{\b_{n-2},n-2}\und z\{\a_n,\a_{n-2}\} =\{\b,\c\}$ we also have
$\und z\{\a_n,\a_{n-2}\}=\{\a_n,\a_{n-2}\}$.  Now
$a_{\b_{n-2},n-2}\{\a_{n-2},\a_n\}=\{\b_{n-2},\a_n\}$ and
$a_{\b,n}\{\b_{n-2},\a_n\} =\{\c,\b\}$. This means
$a_{\b,n}\{\b_{n-2}\}=\{\c\}$.  Also $ a_{\b_{n-2},n}\{\a_n\}=\{\b_{n-2}\}$
and so $a_{\b,n}a_{\b_{n-2},n}\{\a_n\}=\{\c \}$.

Now Corollary~\ref{betaimage} and Lemma~\ref{lm:enrn*} give
$$a_{\b,n}a_{\b_{n-2},n}e_{n}e_{n-2} \homog a_{\c,n}e_n\und z'e_{n-2} \homog
a_{\c,n}e_na_{\c_{n-2},n-2}e_{n-2}\und z$$ for some $\und z \in Z_{Y(2)}$
and $\und z' \in Z_n$.  In particular the lemma holds in this case.

Suppose then that $n$ is in the support of neither $\b$ nor $\c$.
We argue by induction on $n$.  The two reductions of $\und a
e_ne_{n-2}$ are $a_{\b,n-1}e_{n-1}e_na_{\b_{n-2},n-2}e_{n-2}$ and
$a_{\c,n-1}e_{n-1}e_na_{\c_{n-2},n-2}e_{n-2}$ up to right multiples
by elements of $Z_{Y(2)}$.

We know both $\b_{n-2}$ and $\c_{n-2}$ do not have $n$ or $n-1$ in their
support.  We will argue that neither has $n-2$ in the support either.
Then $a_{\b_{n-2},n-2}e_{n-2} \homog a_{\b_{n-2},n-3}e_{n-3}e_{n-2}$ and
similarly for $a_{\c_{n-2},n-2}e_{n-2}$.

Notice $\{\b,\c\}a_{\b,n-1} = \{\c', \a_{n-1}\}$ and $\{\b,\c\}a_{\b,n-1}e_n
= \{\b_{n-2},\a_n\}$.  Now $\c'$ is a root lying in the subsystem of type
$\c' \in D_{n-1}$ and so does not have $n$ in its support.  Moreover,
$(\a_{n-1},\c')=0$ and $\c' \neq \a_{n-1}^* $ as $\{\b,\c\}$ is in the
$W$-orbit of $B_{Y(2)}$. This means $\c'$ does not have $n-2$ in its support
either or $(\c',\a_{n-2})=-1$, in which case, by definition, $\c'=\b_{n-2}$
and this root does not have $n-2$ in its support as claimed.

The two reductions of $\und a
e_ne_{n-2}$ are now $a_{\b,n-1}e_{n-1}e_na_{\b_{n-2},n-3}e_{n-3}e_{n-2}$ and
 \newline $a_{\c,n-1}e_{n-1}e_na_{\c_{n-2},n-3}e_{n-3}e_{n-2}$ up to right multiples
by elements of $Z_{Y(2)}$.

Now both $a_{\b,n-1}e_{n-1}a_{\b_{n-2},n-3}e_{n-3}$ and
$a_{\c,n-1}e_{n-1}a_{\c_{n-2},n-3}e_{n-3}$ belong to $F_{n-1}$.  By
induction on $n$, one can be reduced to the other---up to a right factor
from $Z_{\{\a_{n-1},\a_{n-3}\}}$, say
\begin{eqnarray}
\label{eq:btoc}
a_{\b,n-1}e_{n-1}a_{\b_{n-2},n-3}e_{n-3}&\homog&
a_{\c,n-1}e_{n-1}a_{\c_{n-2},n-3}e_{n-3}\und z'
\end{eqnarray}
for $\und z'\in
Z_{\{\a_{n-1},\a_{n-3}\}}$.  Due to (HNeee) and the definition of $Z_{n-1}$
we have $\und{z}_{n-1}^*e_{n-2}\isog e_{n-1}\und{z}_{n-2}^*$.  Terms in
$\und z'$ generated by $e_i$ or $r_i$ with $i<n-3$ are in $Z_{Y(2)}$.  If
there is $\und{z}_{n-3}^*$, then, using Lemma~\ref{lm:enrn*}(iv), we can
replace it with $\und{z}_{n-1}^*$.  If $\und z'$ is a product of generators
with index less than $n-4$ we have $\und z'e_ne_{n-2}\isog e_ne_{n-2}\und
z'$ with $\und z'\in Z_{Y(2)}$.  In case $\und z'=\und{z}_{n-1}^*$, we find
$\und{z}_{n-1}^*e_ne_{n-2}\homog\und{z}_{n-1}^*e_{n-2}e_{n}\homog
e_{n-1}\und{z}_{n-2}^*e_{n-2}e_n\delta^{-1}\homog
e_{n-1}e_{n-2}e_{n}\delta^{-1}\und{z}_{n-2}^*$,
as $e_{n-2}]\homog e_{n-2}^2\delta^{-1}$.

Now multiplication by
$e_ne_{n-2}$ on both sides of the reduction (\ref{eq:btoc})
and application of
Lemma \ref{lm:enrn*}(iv) gives
\begin{eqnarray*}
a_{\b,n}e_{n}a_{\b_{n-2},n-2}e_{n-2}&\homog&
a_{\c,n-1}e_{n-1}a_{\c_{n-2},n-3}e_{n-3}\und z'e_{n-2}e_n\\
&\homog& a_{\c,n-1}e_{n-1}a_{\c_{n-2},n-3}e_{n-3}e_{n-2}e_n\und z\\
&\homog& a_{\c,n}a_{\c_{n-2},n-2}\und z
\end{eqnarray*}
for $\und z \in Z_n$, as required.  Here, $\und z$ is the same as $\und z'$
unless $\und z'$ has a factor $\und{z}_{n-1}^* $ or $\und{z}_{n-3}^*$, in
which case we can take it to be $\und z_{n-1}^*$ by Lemma~\ref{lm:enrn*}.
In the case with $\und{z}_{n-1}^*$ occurring, the extra $e_{n-1}$
commutes to the left.
If $\und z'=\und z''\und{z}_{n-1}^*$, then $\und z$ reduces to $\und
z''\und{z}_{n-2}^*$.

Next we deal with $\und b$. In this case again, using Proposition
\ref{lm:Theta1}, we can proceed as above. The words
$\und b e_{Y^*(2)}$ can be taken to belong to the
subalgebra $Q$ of type $\A_{n-1}$ and so the reduction is simpler.
Notice here $\{\b,\b^*, \c,\c^*\}a_{\b, n}e_n=\{\c', \c'^*,\a_n,\a_n^*\}$
with support of $\c'$  in $\D_{n-2}$.
\end{proof}

\np This case, for just two roots, extends to admissible sets of arbitrary
size by the next lemma.

\begin{Lm}\label{lm:allthesame}
Let $\und a$ be a word in $F_n$ and choose $Y\in\MY$ such that
$\und a(\emptyset) \in WB_Y$. Let $t\in\{0,\ldots,\lfloor n/2\rfloor\}$ be
such that $Y \in \{Y(t),Y'(t),Y^*(t)\}$.
Then there are positive roots
$\b_{n-2k}$ for $k = 0,\ldots,t-1$ such that $\b_{n-2k}$ has support in
$\D_{n-2k}$ for each $k$ and $\und a \he_Y$ can be reduced to an element of
$\und a' Z_Y$ where
$$\und a' = a_{\b_n,n} a_{\b_{n-2},n-2}  \cdots a_{\b_{n-2t+2},n-2t+2}$$
and every word reduced from $\und a \he_Y$
as in Lemma \ref{lm:mainReduction}
can also be reduced to a word in $\und a' Z_Y$.
\end{Lm}

\nl
\begin{proof}
Set $B = \und a B_Y$.
By Lemma \ref{lm:mainReduction} there is a unique reduction up to right
multiplication by elements of $Z_Y$ for each ordering of the elements of
$B$.  We use Lemma~\ref{lm:twoRoots} to see that the order of, say the first
two, does not matter, in the sense that one reduction can be reduced to
another.  Continuing this way with $\a_{n-2}$ and $\a_{n-4}$, we see that
the words as in Lemma \ref{lm:mainReduction} for all orders of the roots of
$B$ can be reduced to a particular one.  This proves the lemma.
\end{proof}

\begin{Notation}\label{not:ABN} \rm
The lemma allows us to define $a_{B,n}$, for $B\in WB_Y$, as the unique word
$\und a' e_Y\in F_n$ up to homogeneous equivalence and powers of $\delta$
determined by Lemma \ref{lm:allthesame} with $B = \und a' B_Y = \und
a'(\emptyset)$. When $t=0$, we take $a_{B,n}$ to be the identity of the
Brauer algebra.
\end{Notation}

\np The sets $U_Y$ were introduced in Proposition
\ref{prop:ZY}(\ref{not:VYt}).

\begin{Thm}\label{thm:notthetaunique}
Each $\und a\in F_n$ can be reduced to a word of the form $a_{B,n}\und
za_{B',n}^{\op}\delta^k$ where $k$ is an integer, $B$, $B'\in WB_Y$, and
$\und z\in U_Y$ for some $Y\in\MY$.
\end{Thm}

\nl
\begin{proof}
Put $B = \und a(\emptyset)$ and $B' = {\und a}^{\op}(\emptyset)$.  It
follows from Lemma~\ref{lm:mainReduction} that the two sets belong to the
same $W$-orbit inside $\AO$, namely the one containing $B_Y$.  It suffices
to prove the statement of the theorem for $Z_Y$ instead of $U_Y$ because
$\und a (\emptyset) = B$ and, by Proposition \ref{prop:BrMonAction}, the
presence of $e_Ye_i$ in $\und z$ for some $i$ non-adjacent to all members of
$Y$ in $M$ would imply that $a_{B,n}\und za_{B',n}^{\op}(\emptyset) $
contains $a_{B,n}(B_Y\cup\{\a_i\})$, a set of size greater than $|B|$;
however $\und a (\emptyset)=B$ has size $|B|$, a contradiction.

Consider the case $Y = Y(t)$ and suppose $B\in W B_{Y(t)}$.  If $B =
\emptyset$, then $\und a$ does not contain any occurrences of $e_i$ as
$e_i(\emptyset)$ contains $\a_i$ (cf.\ the last assertion of Proposition
\ref{prop:BrMonAction}).  This means that $\und a$ is a product of $r_i$ and
the Matsumoto--Tits rewrite rules for $W$ suffice for the validity of the
theorem in this case, with $t=0$ and $Y = Y(0) = \emptyset$.

Therefore, we may assume that $B\ne\emptyset$, so there is an index $i$ such
that $e_i$ occurs in $\und a$.  If $i\ne n$, then by homogeneous equivalence,
we can replace $e_i$ by $e_{i,n} e_{n-1,i}$.  Thus $\und a = \und b e_n \und
c$ for certain $\und b$, $\und c\in F_n$. By Lemma \ref{lm:allthesame}
applied to both $\und b$ and ${\und c}^{\op}$, we can reduce $\und a$ to
$a_{\b,n} \und{z}_n a_{\b',n}^{\op}$ for some $\b\in B$, $\b'\in B'$ and
$\und{z}_n\in Z_n$.  Then, by an argument as in the proof of
Lemma~\ref{lm:mainReduction}, $\und{z}_n =e_n \und a' (\und{z}_n^*)^\eps$ with
$\und a' \in F_{n-2} $ and $\eps\in\{0,1\}$.
This deals with the case where $t=1$.

Suppose $t>1$.  By induction on $n$, the word $\und a'$ reduces to $
a_{D,n-2} \und z' a_{D',n-2}^{\op}$ for some $\und z'\in Z_{X''}$, where $X'' =
\{\a_{n-2},\ldots,\a_{n-2t+2}\}$ and $D$ and $D'$ are admissible sets in the
root system of type $\D_{n-2}$ with support in $\{1,\ldots,n-2\}$.  Due to
Lemma~\ref{lm:enrn*} (ii), (iv), $$e_ne_{n-2}\und{z}_n^*\isog
e_ne_n\und{z}_{n-2}^*\isog \delta e_n\und{z}_{n-2}^*\isog
e_ne_{n-2}\und{z}_{n-2}^*.$$ By induction on $t$, this gives
$e_{Y(t)}\und{z}_n^*\isog e_{Y(t)}\und{z}_{n-2t+2}^*$, which is the same as
$e_ne_{X''}\und{z}_n^*\isog e_ne_{X''}\und{z}_{n-2t+2}^*$.  So, by Lemmas
\ref{lm:allthesame} and \ref{lm:enrn*} parts (iii) and (iv),
\begin{eqnarray*}
\und a&\homog& a_{\b,n} \und{z}_n a_{\b',n}^{\op} \isog a_{\b,n}e_n \und a'
(\und{z}_n^*)^\eps a_{\b',n}^{\op} \homog a_{\b,n}e_n a_{D,n-2} \und z'
a_{D',n-2}^{\op} (\und{z}_n^*)^\eps a_{\b',n}^{\op}\\ &\homog&a_{\b,n}
a_{D,n-2} e_n\und z'(\und{z}_n^*)^\eps a_{D',n-2}^{\op} a_{\b',n}^{\op} \\
&\homog& a_{\b,n}a_{D,n-2} e_n \und z' (\und{z}_{n-2t+2}^*)^\eps
a_{D',n-2}^{\op} a_{\b',n}^{\op} \homog a_{B,n} z a_{B',n}^{\op}
\end{eqnarray*}
for some $\und z\in Z_{Y(t)}$. This handles the case $Y = Y(t)$.

If $B\in WB_{Y'(n/2)}$, the same arguments apply. Finally, if $B\in
WB_{Y^*(t)}$, then, due to
Proposition \ref{lm:Theta1},
the same arguments apply to $Q$ with the root system
of type $\A_{n-1}$ having support in $\{2,\ldots,n\}$.
\end{proof}

\np The following corollary extends Lemma \ref{lm:Matsumoto} to the
the Brauer monoid.

\begin{Cor}\label{cor:finalreduction}
For each $\und a\in F_n$, all reduced elements of $F_n$ reducible from $\und
a$ are homogeneously equivalent to an element of the form $a_{B,n}\und z
a_{B',n}^{\op}\delta^k$ with $B$ and $B'$ in $W B_Y$ for some $Y\in\MY$ and
$\und z\in U_Y$. Here the elements $B$ and $B'$ are uniquely determined by $B =
\und a(\emptyset)$ and $B' = {\und a}^\op(\emptyset)$, respectively.
\end{Cor}

\nl
\begin{proof}
The form is immediate from the theorem. Uniqueness up to homogeneous
equivalence follows from Lemma \ref{lm:allthesame} for $a_{B,n}$ and
$a_{B',n}$ and from the Matsumoto--Tits' rewrite rules for $\und z\in
U_Y$, as stated in Proposition \ref{prop:ZY}.
\end{proof}

\np
As a consequence of Corollary~\ref{cor:finalreduction},
all reduced words that are reductions from $\und a \in F_n$ in
$\Br(\D_n)$ are homogeneously equivalent.

The proof of Corollary~\ref{cor:finalreduction} has implications for the
ordinary Brauer algebra of type $\A_{n-1}$ which we can take to be generated
by $r_i$, $e_i$ for $2\leq i \leq n$.  Here there are no $r_i^*$, and
$U_{Y(t)}$ consists of the reduced words on $\und{s}_i$ for $2\le i\leq
n-2t$, while $W$ is generated by the $r_i$ for $2\le i \le n $ and is
isomorphic to the symmetric group of $n$ points.

\begin{Cor}\label{cor:finalreductionAn}
Let $\Br(\A_{n-1})$ be the Brauer algebra of type $\A_{n-1}$.  Let $\pi$ be
the map from $F_{n-1}$ to $\Br(\A_{n-1})$ taking $r_i$ or $e_i$ to the
element in $\Br(\A_{n-1})$ with the same label.  For each $\und a\in
F_{n-1}$, all reduced words in $F_n$ reducible from $\und a$ are
homogeneously equivalent to an element of the form $a_{B,n-1}\und
za_{B',n-1}^{\op}\delta^k$ with $B$ and $B'$ in $W B_{Y(t)}$ for some
$t\in\{0,1,\ldots,\lfloor \frac{n-1}{2}\rfloor\}$ and $\und z\in
U_{Y(t)}$. Here the elements $B$ and $B'$ are uniquely determined by $B =
\und a(\emptyset)$ and $B' = {\und a}^\op(\emptyset)$, respectively.  Also,
$U_{Y(t)}$ is the Weyl group of type $\A_{n-1-2t}$.
\end{Cor}

\np
We chose the index $n-1$ here so there will be no confusion between these
coefficients and the ones in Corollary~\ref{cor:finalreduction}.

\section{Proof of Theorem~\ref{th:main} and Corollary~\ref{TempLieb}}
\label{sec:conclusion}
In this section we prove Theorem \ref{th:main} and Corollary~\ref{TempLieb}.

\np {\bf Proof of Theorem~\ref{th:main}.}
By Lemma~\ref{lm:basis}, the rank of $\BMW(\D_n)$
is at least
$\dim(\Br(\D_n))$, which by \cite[Theorem 1.1]{CFW} equals $(2^n+1)
n!!-(2^{n-1}+1)n!$.

\np Now let $T$ be the set of elements $a_{B,n}\und za_{B',n}^{\op}$ in
$F_n$ as in Corollary~\ref{cor:finalreduction}.  Then the elements of $T$
correspond to triples $(B,B',\und z)$ where $B$ and $B'$
are in the $W$-orbit in $\AO$ containing $B_Y$ for some $Y\in \MY$ and $\und
z\in U_Y$.  Now $T$ is a finite set and, by Corollary
\ref{cor:finalreduction}, every $\und a \in F_n$ reduces to an element of
$T$ up to a power of $\delta$.

For the remainder of the proof of the first statement of Theorem
\ref{th:main}, we note that, by \cite[Proposition 4.9 and the proof of
Theorem 1.1]{CFW}, $(\pi(t))_{t\in T}$ is a basis of $\Br(\D_n)$.  Now
Proposition \ref{prop:BMWbasis} applies, so $(\rho(t))_{t\in T}$ is a basis of
$\BMW(\D_n)$.  This shows that $\BMW(\D_n)$ is free of rank as claimed in
Theorem~\ref{th:main}.

To show that $\BMW(\D_n)$ tensored over $\Q(l,\delta)$ is semisimple we use
the surjective equivariant map $\mu\,\colon\BMW(\D_n)\otimes _R
\Q(\delta)[l^{\pm1}]\to\Br(\D_n)$ over $\Q(\delta)$; cf.~Definitions
\ref{df:BrMonoid}.  We know its image $\Br(\D_n)$ is semisimple by
\cite[Corollary~5.6]{CFW} and so has no nilpotent left ideals.  Suppose
$\BMW(\D_n)\otimes _R\Q(\delta,l)$ has a nontrivial nilpotent ideal.  Take a
nonzero element of it expressed in the basis we have found.  Multiply the
element by a suitable polynomial in $l$ so that all coefficients are in
$\Q(\delta)[l^{\pm 1}]$.  As in the proof of Lemma~\ref{lm:basis}, rescale
the coefficients by a power of $l-1$ so that all coefficients remain in
$\Q(\delta)[l^{\pm1}]$ but some coefficient $\lambda_s$ lies outside
$(l-1)\Q(\delta)[l^{\pm1}]$.  The result is a nonzero nilpotent element in
$\BMW(\D_n)\otimes \Q(\delta)[l^{\pm1}]$ with $\mu(\lambda_s) \neq 0$, so
its image under $\pi$ is a nonzero nilpotent element of $\Br(\D_n)$.
Furthermore, any multiple is nilpotent both in $\BMW(\D_n)\otimes
\Q(\delta,l)$ and in $\Br(\D_n)$ and so generates a nontrivial nilpotent
ideal of $\Br(\D_n)$, a contradiction.  This completes the proof of Theorem
\ref{th:main}.

\np Although we did not need the statement for the proof of the main
theorem, it may be worthy of mention that, by Proposition
\ref{prop:BMWbasis}, each word in $T$ as above is reduced.

We will need the elements $\rho(a_{B,n})$ for the words $a_{B,n}$ in
$F_n$ introduced in Notation \ref{not:ABN}. These words were defined up
to homogeneous equivalence. Since different elements from the homogeneous
class of $a_{B,n}$ may give different elements in $\BMW(\D_n)$, see Remark
\ref{rmk:notInBMW}, we need to select a particular element in each class.

\begin{Notation}\label{not:rBn} \rm
Let $Y\in\MY$.
For each $B\in WB_Y$, we take $a_{B,n}$ to be a specific word in $F_n$ from
its homogeneous equivalence class in $F_n$ and write
$\rr_{B,n}=\rho(a_{B,n})$ for its image in $\BMW(\D_n)$ under $\rho$.
\end{Notation}

\begin{Cor}\label{cor:BDbasis}
For $n\ge4$, the elements $\rr_{B,n}\rho(\und z)\rr_{B',n}$ for $(B,B',\und
z)\in \bigcup_{Y\in\MY} WB_Y \times WB_Y \times U_Y$ are a basis of
$\BMW(\D_n)$.
\end{Cor}

\begin{Remark}\label{Theta'}
\rm
Let $\Theta'$ be the ideal of
$\BMW(\D_n)$ generated by $e_1e_2$. We can
also choose a basis of $\Theta'$ of the form $\rr_{B,n-1}\rho(\und
z)\rr_{B',n-1}$ where the $\rr_{B,n-1}$ are chosen as in
Notation~\ref{not:rBn} using Corollary~\ref{cor:finalreductionAn} for
$\Br(\A_{n-1})$.
\end{Remark}

\begin{Remark}\label{substructures}
\rm A consequence of Theorem \ref{th:main} is that natural subalgebras
generated by $\{g_i,e_i\mid i\in K\}$ for $K$ a set of nodes of $M$ have the
usual desired subalgebra structure, that is, are naturally isomorphic to the
BMW algebra whose type is the restriction of $M$ to $K$.  In particular, the
subalgebra generated by $\{ g_i,e_i\mid 2\leq i\leq n\} $ is the full
$\BMW(\A_{n-1})$ rather than a proper homomorphic image.  The same applies to
the algebra generated by all $g_i$, $e_i$ for $i\leq n-1$ which is
$\BMW(\D_{n-1})$ and not a proper image.
\end{Remark}

\medskip {\bf Proof of Corollary~\ref{TempLieb}.}
The generalized Temperley--Lieb algebra of type
$\D_n$ has been studied in \cite{Fan,Gra,Gre}.  The elements $e_i$
either in $\BMW(\D_n)$ or in $\Br(\D_n)$ commute for $i\not \sim j$ by (HCee).
For $i\sim j$, we have $e_ie_je_i=e_i$ by (HNeee).  Also, $e_i^2=\delta e_i$
by (HSee).  The free algebra on $e_1,\ldots,e_n$ with this presentation over
$\Z[\delta^{\pm1}]$ is called the (generalized) Temperley--Lieb algebra of
type $\D_n$ over $\Z[\delta^{\pm1}]$; we will denote it by $\TL(\D_n)$.  The subalgebra generated by
$e_1,\ldots,e_n$ in $\Br(\D_n)$ is a homomorphic image of $\TL(\D_n)$; the
subalgebra of $\BMW(\D_n)$ generated by these elements is a homomorphic image
of $\TL(\D_n)\otimes_{\Z[\delta^{\pm1}]}R$.  The words in $F_n$ corresponding
to generators for these subalgebras consist solely of the symbols
$e_1,\ldots,e_n,\delta$ and so are of height~$0$.

In \cite[Theorem~4.2 and Lemma~6.5]{Gre}, a description of a generating set
for the Temperley--Lieb algebra is given in terms of decorated diagrams with
some restrictions.  In \cite{CGW3} diagrams such as these were introduced
for the full algebra $\Br(\D_n)$.  In particular, in
\cite[Theorem~1.1]{CGW3} it is shown there is an isomorphism, $\nu$, from
$\Br(\D_n)$ to the span of the diagrams as a basis over
$Z[\delta,\delta^{-1}]$.  In \cite[Lemmas~6.5 and 6.6]{Gre} it is shown that
the specific images $\nu(e_i)$ generate the full Temperley--Lieb algebra
and so the $e_i$ in $\Br_n(\D_n)$ generate the full Temperley--Lieb algebra.
The actual multiplication of $\nu(e_i)$ with $\nu(e_j)$ in
\cite[Section~4]{CGW3} has a coefficient $\xi$ which sometimes appears.
However, by results in \cite{CGW3} the coefficient $\xi$ does not appear for
words of height~$0$ and so does not appear here.

Now apply Proposition~\ref{lm:basis} to see that the algebra generated by
$\rho(e_i)$ is the full Temperley--Lieb algebra over $R$.

This completes the proof of Corollary~\ref{TempLieb}.  It follows from this
that the subalgebra of $\Br(\D_n)$ generated by $e_j$ $(j\ge2)$ is
isomorphic to the Temperley--Lieb algebra of type $\A_{n-1}$.

In \cite{CW4}, it is shown that the Temperley--Lieb monomials are the terms
$a_{B,n}$ for $B$ of height $0$, where the concept of height $0$ for $B$ is
given in Section~\ref{sec:discussion}.

\begin{Remark}
\label{rmk:reps}
\rm
By use of $\mu$ and the Tits Deformation Theorem, see \cite[IV.2, exercise
26]{Bour} or \cite[Lemma 85]{Stein}, it can be shown that the irreducible
degrees associated to $\BMW(\D_n)$ are the same as for $\Br(\D_n)$.
This can
also be shown by use of Theorem~\ref{thm:notthetaunique} for representations
with $\Theta' $ in the kernel as in \cite{CFW} and for the others from the
connection of $\Theta'$ to $\BMW(\A_{n-1})$ as in the proof of
Theorem~\ref{th:main}.
\end{Remark}

\section{Cellularity}\label{cellular}
Let $S$ be a commutative algebra over $R$.  In this section we prove
Theorem~\ref{th:cellular}, which states that $\BMW(\D_n)\otimes_R S$ is
cellular in the sense of Graham--Lehrer \cite[Definition~1.1]{GL} if
$S$ contains an inverse to $2$.
We recall the definition from \cite{GL}.

\begin{Def}
\label{df:cellular}
\rm
An associative
algebra $\alg$ over a commutative ring $S$ is \emph{cellular}
if there is a quadruple
$(\Lambda, T, C, *)$ satisfying the following three conditions.

\begin{itemize}
\item[(C1)] $\Lambda$ is a finite partially ordered set.  Associated to each $\lambda \in \Lambda$, there is
a finite set $T(\lambda)$.  Also, $C$ is an injective map
$$ C:\ \ \ \coprod_{\lambda\in \Lambda} T(\lambda)\times T(\lambda)
\rightarrow \alg$$ whose image is an $S$-basis of $\alg$.

\item[(C2)]
The map $*:\alg\rightarrow \alg$ is an
$S$-linear anti-involution such that
$C(x,y)^*=C(y,x)$ whenever $x,y\in
T(\lambda)$ for some $\lambda\in \Lambda$.

\item[(C3)] If $\lambda \in \Lambda$ and $x,y\in T(\lambda)$, then, for any
element $a\in \alg$,
$$aC(x,y) \equiv \sum_{u\in T(\lambda)} r_a(u,x)C(u,y) \
\ \ {\rm mod} \ \alg_{<\lambda},$$ where $r_a(u,x)\in S$ is independent of $y$
and where $\alg_{<\lambda}$ is the $S$-submodule of $\alg$ spanned by $\{
C(x',y')\mid x',y'\in T(\mu)\mbox{ for } \mu <\lambda\}$.
\end{itemize}
Such a quadruple $(\Lambda, T, C, *)$ is called a {\em cell datum} for
$\alg$.
\end{Def}

\np
Now let $S$ be an integral domain containing $R$ as a subring with
$2^{-1}\in S$.  We introduce a quadruple $(\Lambda, T,C,*)$ and prove that
it is a cell datum for $A = \BMW(\D_n)\otimes_R S$.  The map $*$ on $A$ will
be the opposition map ${\cdot}^{\op}$ of Notation \ref{not:op}.  Before
describing the other three components of the quadruple $(\Lambda, T, C, *)$,
we will relate the subalgebras of $A$ generated by monomials corresponding
to the elements of $U_Y$, defined in Proposition \ref{prop:ZY}(iii), to
Hecke algebras.  For this purpose we need a version of
Proposition~\ref{prop:ZY} that applies to $A$ rather than $\BrM(\D_n)$. This
requires a version of Lemma \ref{lm:enrn*} for $\BMW(\D_n)$ rather than
$F_n/\isog$, with $\homog$ replaced by equality in $\BMW(\D_n)$.  Here, as
in Remark~\ref{Theta'}, we let $\Theta'$ be the ideal of $\BMW(\D_n)$
generated by $e_1e_2$.

\begin{Lm}\label{lm:engn*} For $n\ge3$, the monomials $\hat{z}_n^* =
\rho(\und{z}_n^*)$ in $\BMW(\D_n)$
satisfy the following
equations, where $g_n^* = \rho(r_n^*)$.  Here we are considering $e_i$ and
$g_i$ to be in $\BMW(\D_n)$, namely $g_i=\rho(r_i)$ and $e_i=\rho(e_i)$
where $r_i\in F_n$ and the $e_i$ within the parentheses is also in $F_n$.
\begin{enumerate}[(i)]
\item
$g_n^*e_n = \hat{z}_n^* = e_ng_n^* $.
\item
$\hat{z}_n^* = e_{n,3}g_2e_1e_{3,n}$.
\item
\label{lm:xinB}
For $n\ge4$ and
$i\in\{1,\ldots,n-2\}$, $e_i\hat{z}_n^* =  \hat{z}_i^* e_n = e_n\hat{z}_i^*$
and $g_i\hat z_n^* =  \hat{z}_n^*g_i$.
\item
$\hat{z}_n^*e_{n-2}
= e_n\hat{z}_{n-2}^*$ and $e_n\hat{z}_n^*= \hat{z}_n^*e_n= \delta \hat{z}_n^*$.
\item
\label{en*sqaredB}
$(\hat{z}_n^*)^2 \in \delta e_n - m \delta \hat{z}_n^*+\Theta'$.
\end{enumerate}
For $n$ equal to $1$ or $2$, both (i) and (v) hold.
\end{Lm}

\nl\begin{proof}
Many of the proofs are the same as for Lemma~\ref{lm:enrn*}.  Differences
occur when the relations are not monomial as extra terms occur with
coefficients divisible by $m$.

\nl(i).
The proof is similar to the one of Lemma \ref{lm:enrn*}(i);
note that the relations (RNrre) for $\BMW(\D_n)$ are also binomial.

\nl(ii).  Again, the only relation used in the proof of
Lemma \ref{lm:enrn*}(ii) is (HTeere), which is binomial for $\BMW(\D_n)$.

\nl(iii).
Let $i\in\{2,\ldots,n-2\}$.
The relation $e_i\hat{z}_n^* =
e_n\hat{z}_i^*$ can be derived from
the definition of $e_{k,n}$, and the binomial relations (HCee) and
(HNeee), as in the proof of Lemma \ref{lm:enrn*}.

The proof of $g_i\hat{z}_n^* = \hat{z}_n^* g_i $ is a bit more involved.  By
(HCer), (RNrre), (RSrr) and (RNerr),
\begin{eqnarray*}
g_i\hat{z}_n^* &= &  e_{n,i+2} g_ie_{i+1,2}g_1e_{3,n}
=   e_{n,i+2}g_i^2g_{i+1}e_{i,2}g_1e_{3,n} \\
& = &  e_{n,i+2}g_{i+1}e_{i,2}g_1e_{3,n}
-m e_{n,i+2}g_ig_{i+1}e_{i,2}g_1e_{3,n}\\
&&\qquad\qquad\qquad
+ml^{-1} e_{n,i+2}e_ig_{i+1}e_{i,2}g_1e_{3,n}\\
&=&  e_{n,i+2}e_{i+1}g_{i+2}^{-1}e_{i,2}g_1e_{3,n}
-m e_{n,2}g_1e_{3,n}
+m e_{n,i+2}e_{i,2}g_1e_{3,n}\\
&=&  e_{n,i+1}g_{i+2}^{-1}e_{i,2}g_1e_{3,n}
-m \hat{z}_n^*
+m e_{i,2}g_1e_{n,i+2}e_{3,n}\\
&=&  e_{n,2}g_{i+2}^{-1}g_1e_{3,n}
-m \hat{z}_n^*
+m e_{i,2}g_1e_{n,i}e_{n}\\
&=& e_{n,2}g_{i+2}^{-1}g_1e_{3,n}
-m \hat{z}_n^*
+m \hat{z}_{i}^*e_{n}.
\end{eqnarray*}
Since each of the three summands is invariant under opposition,
(observe that $(\hat{z}_n^*)^\op = \hat{z}_n^* $ follows from (i)),
so is $g_i\hat{z}_n^*$. This shows
$g_i\hat{z}_n^* = (g_i\hat{z}_n^*)^\op = (\hat{z}_n^*)^\op g_i
= \hat{z}_n^* g_i$.

The case $i=1$ is notationally different but can be done the same way
as $i=2$.

\nl(iv).
By (iii) with $i=n-2$ we have
$e_{n-2}\hat{z}_n^* =  \hat{z}_{n-2}^* e_n = e_n\hat{z}_{n-2}^*$.
Taking images under $\cdot^\op$ and using opposition invariance  of
$\hat{z}_n^*$, we find
$\hat{z}_n^* e_{n-2}
=(e_{n-2}\hat{z}_n^* )^\op
=(\hat{z}_{n-2}^*e_n)^\op
= e_n\hat{z}_{n-2}^*$, as required for the first equation.
The second chain of equations is a direct consequence of (RSee).

\nl(v).

For $n\ge3$, by (HSee), (HCer), (HNeee), and (RSrr),
\begin{eqnarray*}
(\hat{z}_n^*)^2 &=&
e_{n,3}g_1e_2e_{3,n}e_{n,3}g_1e_2e_{3,n} = e_{n,3}g_1e_2e_3g_1e_2e_{3,n}
\delta \\
& = & \delta  e_{n,3}g_1e_2e_3e_2g_1e_{3,n}
= \delta e_{n,3}g_1e_2g_1e_{3,n}
= \delta e_{n,3}g_1^2e_2e_{3,n}
\\
&=&  \delta  e_{n,3}(1 -mg_1+ml^{-1}e_1)e_{2,n}\\
&=&  \delta  e_{n,3}e_{2,n} -\delta me_{n,3}g_1e_{2,n} + \delta ml^{-1}e_{n,3}e_1e_2e_{3,n}
\\
&\in&  \delta e_n -   \delta m\hat{z}_n^*  +\Theta'.
\end{eqnarray*}

 The cases $n=1$ and $n=2$ are easily proved.
\end{proof}

\begin{Def}\label{df:Jt} \rm
For $t\in \{0,\ldots,\lfloor n/2\rfloor\}$, we write $J_{t+1}$ to denote the
ideal of $A$ generated by $\Theta'$, $e_{Y(t')}$ for all $t'>t$, and
$e_{Y'(n/2)}$ if $n$ is even and $n>2t$.  In particular, $J_1$ is the ideal
generated by all $e_i$ and, for $n$ even with $t={n}/{2}$, the ideal $
J_{t+1}$ coincides with $\Theta'$.
\end{Def}

\np
Recall the words $\und{s}_i$ $(0\le i \le n-2t)$ and
$\he_{Y(t)}=e_{Y(t)}\delta^{-t}$ given in Notation~\ref{not:zn}.

\begin{Prop}\label{prop:hecke}
Let $t\in\{0,\ldots,\lfloor n/2\rfloor\}$. The monomials
$\hat{s}_i= \rho(\und s_i)$, for $i$ a node of $M_{Y(t)}$,
satisfy the following relations.
\begin{enumerate}[(i)]
\item
The element $\rho(\he_{Y(t)})$ acts as an identity element on the $\hat{s}_i$, that
is, $\rho(\he_{Y(t)})\hat{s}_i = \hat{s}_i $ and $\hat{s}_i\rho(\he_{Y(t)})=\hat{s}_i$,
while $\rho(\he_{Y(t)})^2 =\rho(\he_{Y(t)})$.  Moreover, the $\hat{s}_i$ satisfy the braid
relations (HCrr) and (HNrrr) of Table \ref{BMWTable} with $g_i$ replaced by
$\hat{s_i}$ and $1$ by $\rho(\he_{Y(t)})$.
\item
Each monomial $\hat{s}_i$
satisfies the quadratic Hecke algebra relation
modulo the ideal
$J_{t+1}$, that is,
$\hat{s}_i^2 +m \hat{s}_i - \rho(\he_{Y(t)}) \in J_{t+1}$.
\end{enumerate}
If $n$ is even and $t={n}/{2}$, the corresponding statement holds for
$Y'({n}/{2})$ replacing $Y(n/2)$.
\end{Prop}

\nl
\begin{proof}
Here again $g_i$ and $e_i$ are considered in $\BMW(\D_n)$.

\nl(i).
The relations involving $\rho(\he_{Y(t)})$ are easily derived from Lemma
\ref{lm:engn*}. Note the resemblance with the proof of Proposition
\ref{prop:ZY}.

Use of (RSrr), (HCrr), and
(HNrrr) gives the relations not involving $\hat{s}_0$.
It remains to verify the commuting of $\hat{s}_0$ with
$\hat{s}_i$ for $i\in\{1,\ldots,n-2t\}$.
By Lemma \ref{lm:engn*}(iii)
$g_i\hat{z}_n^* =\hat{z}_n^*  g_i$.
This gives
\begin{eqnarray*}
\hat{s}_{0}\hat{s}_{i}&=&
z_n^*\delta^{-1}\he_{Y(t)}g_i\he_{Y(t)} =
z_n^*g_i\delta^{-1}\he_{Y(t)} =
g_i\he_{Y(t)}z_n^*\delta^{-1} \he_{Y(t)}=
\hat{s}_{i}\hat{s}_{0}.\\
\end{eqnarray*}

\nl(ii).
For $i\in\{1,\ldots, n-2t\}$, we have $\hat{s}_i^2 = g_i\he_{Y(t)}g_i\he_{Y(t)}=g_i^2\he_{Y(t)}
=(1-mg_i+ml^{-1}e_i)\he_{Y(t)} =\he_{Y(t)} - m\hat{s}_i +ml^{-1}
e_{Y(t)\cup\{i\}}\delta^{-1-t}$.  Here $e_{Y(t)\cup\{i\}}$ is in $J_{t+1}$, so
$\hat{s}_i^2 +m\hat{s}_i-\he_{Y(t)}\in
J_{t+1}$.

As for $\hat{s}_0^2$, by Lemma \ref{lm:engn*}(iv), $\hat{z}_n^*e_{Y(t)} =
e_{Y(t)} \hat{z}_{n-2t+2}^*$ and $e_{Y(t)}\hat{z}_n^* = \hat{z}_{n-2t+2}^*
e_{Y(t)} $, so, by Lemma \ref{lm:engn*}(iii),(iv),(v), and in view of
$\Theta'\subseteq J_{t+1}$, we have
\begin{eqnarray*}
\hat{s}_0^2 &=&
\hat{z}_n^*e_{Y(t)} \hat{z}_n^*e_{Y(t)}\delta^{-2t-2}
=
\hat{z}_n^*e_{Y(t)}e_{Y(t)} \hat{z}_n^*e_{Y(t)}\delta^{-3t-2}\\
&=&
e_{Y(t)}(\hat{z}_{n-2t+2}^*)^2e_{Y(t)}^2\delta^{-3t-2}
\\
&\in&
e_{Y(t)}(e_{n-2t+2}- m\hat{z}_{n-2t+2}^*)e_{Y(t)}\delta^{-2t-1}+J_{t+1}
\\
&=&
e_{Y(t)}e_{n-2t+2}e_{Y(t)}\delta^{-2t-1}
- me_{Y(t)} \hat{z}_{n-2t+2}^*e_{Y(t)}\delta^{-2t-1} + J_{t+1}
\\
&=&
\rho(\he_{Y(t)})
- m  \hat{z}_{n}^*e_{Y(t)}e_{Y(t)}\delta^{-2t-1}  + J_{t+1}
\\
&=&
\rho(\he_{Y(t)})
- m  \hat{z}_{n}^*e_{Y(t)}\delta^{-t-1} + J_{t+1}
\\
&=&
\rho(\he_{Y(t)})
- m  \hat{s}_{0}  + J_{t+1}.
\end{eqnarray*}
\end{proof}

\np
We will next exploit the elements $\rr_{B,n}$ of Notation \ref{not:rBn}.
Recall from Proposition \ref{prop:ZY}(iii) the definition of $U_{Y(t)}$.

\begin{Notation}
\rm Let $H_Y$ be the linear span of $\rho(U_Y)$.
\end{Notation}

\np
\begin{Cor}\label{cor:HYt}
For $t\in\{0,\ldots,\lfloor n/2\rfloor\}$, the linear subspace
$H_{Y(t)}$ of $A$ satisfies the following properties.
\begin{enumerate}[(i)]
\item
The linear subspace $H_{Y(t)}+J_{t+1}$ is a subalgebra of $A$
whose quotient algebra mod $J_{t+1}$ is isomorphic to the
Hecke algebra of type $M_{Y(t)}$.  Moreover, $\rho(U_{Y(t)})$
is a basis of $H_{Y(t)}$.
\item For each $i\in\{1,\ldots,n\}$ and $B\in WB_{Y(t)}$, we have
$g_i\rr_{B,n}e_{Y(t)} \in \rr_{r_iB,n}H_{Y(t)}+J_{t+1}$ and
$e_i\rr_{B,n}e_{Y(t)} \in \rr_{B'',n}H_{Y(t)}+J_{t+1}$ for some $B''\in
WB_{Y(t)}$.
\item
\label{opp-invH}
The linear subspace $H_{Y(t)}$ is invariant under opposition.
\end{enumerate}
If $n$ is even, the similarly defined linear span $H_{Y'(n/2)}$ equals
$S\he_{Y'(n/2)}$ and satisfies the same properties.
\end{Cor}

\nl
\begin{proof}
If $\und{a}\in U_{Y(t)} $ is a minimal expression in the $\und{s}_i$ of the
element $\pi(\und{a})\in W(M_{Y(t)})$, then, as a consequence of Lemma
\ref{lm:Matsumoto} and the relations established in Proposition
\ref{prop:hecke}(i), $\rho(\und{a})$ depends only on $\pi(\und{a})$ and not on
the choice of the minimal expression.

\nl(i). By the above and Proposition \ref{prop:hecke}(ii), the spanning set
$\rho(U_{Y(t)})$ of $H_{Y(t)}$ has size at most $|W(M_{Y(t)})|$. Due to
Corollary \ref{cor:BDbasis} there is no collapse, so the spanning set has size
equal to $|W(M_{Y(t)})|$ and is a basis of $H_{Y(t)}$.  By Proposition
\ref{prop:hecke}, the linear subspace $H_{Y(t)}+J_{t+1}$ is closed under
multiplication and satisfies the Hecke algebra defining relations mod
$J_{t+1}$ on the generators $\hat{s}_i$ $(0\le i\le n-2t)$.  In particular,
$(H_{Y(t)}+J_{t+1})/J_{t+1}$ is a quotient of the Hecke algebra of type
$M_{Y(t)}$. But, its rank is equal to $|W(M_{Y(t)})|$, which is the Hecke
algebra dimension, and so $(H_{Y(t)}+J_{t+1})/J_{t+1}$ is isomorphic to the
Hecke algebra of type $M_{Y(t)}$.

\nl(ii). In view of Corollary \ref{cor:BDbasis} and Proposition
\ref{prop:BrMonAction}, $\rr_{B,n}H_{Y,t}\rr_{B',n}+J_{t+1}$ is the
linear span of $J_{t+1}$ and all monomials $x$ in $\BMW(\D_n)$ such that
$\mu(x)(\emptyset) = B$ and $\mu(x)^\op (\emptyset)= B'$.  But $x =
g_i\rr_{B,n}e_{Y,t}$ satisfies $\mu(x)(\emptyset) = r_iB$ and $\mu(x)^\op
(\emptyset) = Y(t)$, so $g_i\rr_{B,n}e_{Y,t}\in
\rr_{r_iB,n}H_{Y,t}+J_{t+1}$.

Similarly, $\mu(e_i\rr_{B,n}e_{Y(t)})(\emptyset) =
\mu(e_i\rr_{B,n})B_{Y(t)} = e_iB $ always contains a member $B''$, say, of
$WB$, and $\mu(e_i\rr_{B,n}e_{Y(t)})^\op (\emptyset) =
e_{Y(t)}(\pi(a_{B,n}))^\op\{\a_i\} $ contains $B_{Y(t)}$, so
$e_i\rr_{B,n}e_{Y(t)}\in \rr_{B'',n}H_{Y(t)}+J_{t+1}$.  Here
if $\a_i \perp B$ the expression is in $J_{t+1}$.

\nl(iii).  It is readily verified that each $\hat{s}_i$ is fixed under
opposition.  As the opposite of a minimal expression in the $\hat{s}_i$ is
again a minimal expression, $\rho(U_{Y(t)})$ is invariant under opposition.
Hence, so is $H_{Y(t)}$.
\end{proof}

\np We now give the cell datum for $A = \BMW(\D_n)_R\otimes S$.  View
$\A_{n-1}$ as the subdiagram of $\D_n$ on the nodes $2,\ldots,n$.  As an
algebra over $S$, the ideal $\Theta'$ of $A$ generated by $e_1e_2$ is
isomorphic to the ideal of $\BMW(\A_{n-1})\otimes_R S$ generated by $e_2$;
see Proposition~\ref{lm:Theta1}.  The ideal generated by $e_2$ is a cellular
algebra as $\BMW(\A_{n-1})$ is cellular by \cite[Theorem~3.11]{Xi} and it
inherits the cellular structure from that of $\BMW(\A_{n-1})$.  In fact, it
corresponds to the ideals with cell datum associated with partitions of
$n-2t$ for $1\le t \le \lfloor n/2 \rfloor$.

Let
$(\Lambda_\theta,T_\theta,C_\theta,*_\theta)$ be the cell datum for
$\Theta'$. It is clear from \cite[Theorem~3.11]{Xi} that $*_\theta$
coincides with the restriction to $\Theta'$ of the map ${\cdot}^\op$.
Moreover, the elements $g_1-g_2$ and $e_1-e_2$ are in the kernel of
the action of $A$ on $\Theta'$ by left multiplication,
as well as by right multiplication.

For $0\leq t \leq \lfloor n/2 \rfloor$ we let $(\Lambda_t,T_t,C_t,*_t)$ be the
cell datum for the Hecke algebra $H_{Y(t)}$ mod $J_{t+1}$ of type
$M_{Y(t)}$ (see Corollary \ref{cor:HYt}(i))
with $*_t$ the restriction to $H_{Y(t)}$ of ${\cdot}^\op$. If
$n=2t$, there is another copy needed which we denote $(\Lambda'_{n/2},
T'_{n/2},C'_{n/2},*'_{n/2})$; it corresponds to the admissible set $Y'(n/2)$.
By \cite{MG}, these cell data are known to exist if
$\frac{1}{2}\in S$.
We take the values of $C_t$ to be in $H_{Y(t)}$.

The poset $\Lambda$ is the disjoint union of $\Lambda_\theta$ together with
the posets $\Lambda_t$ of the
cell data for the various Hecke algebras $H_{Y(t)}$ mod $J_{t+1}$, as well as
$\Lambda'_{n/2}$ if $n$ is even.  We make $\Lambda$ into a poset as follows.
For a fixed $t$ or $\theta$ it is already a poset, and we keep the same
partial order.  Furthermore, any element of $\Lambda_t$ is greater than any
element of $\Lambda_s$ if $t<s$.  In particular the elements of $\Lambda_0$
are greater than the elements of $\Lambda_t$ for any $t\ge 1$.  Moreover, if
$n$ is even, any element of $\Lambda'_{n/2}$ is smaller than any element of
$\Lambda_{t}$ for $t<n/2$.  Finally, we decree that any element of
$\Lambda_\theta$ is smaller than any element of $\Lambda_{t}$ $(0\le t\le
n/2)$ or $\Lambda'_{n/2}$.

\np Let $ t\in\{0,\ldots,\lfloor n/2 \rfloor\}$.  For
$\lambda \in \Lambda_t$, we set $T(\lambda)= WB_{Y(t)} \times T_t(\lambda)$ and,
if $n$ is even, for $\lambda \in \Lambda'_{n/2}$, we set $T(\lambda)= WB_{Y'(n/2)}
\times T'_{n/2}(\lambda)$.
For $\lambda \in \Lambda_\theta$, we set $T(\lambda)=
T_\theta(\lambda)$.
This determines $T$.

\np
We define $C$ as follows. For $t\in\{0,\ldots,\lfloor n/2\rfloor\}$,
$\lambda\in\Lambda_t$,  and  $(B,x),(B',y)\in T(\lambda)$, we have
$$C\big((B,x),(B',y)\big)=\rr_{B,n}C_t(x,y)\rr_{B',n}^\op.$$
Similarly on $\Lambda'_{n/2}\times\Lambda'_{n/2}$.
For $\lambda\in\Lambda_\theta$, the map
$C$ on $T(\lambda)\times T(\lambda)$
is just $C_\theta$.

\np
Since we have already defined $*$ by the opposition map,
this concludes the definition of $(\Lambda, T, C, *)$.
We next verify the conditions (C1), (C2), (C3).

\nl(C1).  The map $C$ has been chosen so that its image is an $S$-basis of
$\Theta'$ (the image of $C_\theta$), joint with the set of elements
$\rr_{B,n}C_t(x,y)\rr_{B',n}^\op$ for $B,B'\in WB_{Y(t)}$ and $C_t(x,y)$
running through a basis of $H_{Y(t)}$, and
$\rr_{B,n}C'_{n/2}(x,y)\rr_{B',n}^\op$ for $B,B'\in WB_{Y'(n/2)}$ and
$C'_{n/2}(x,y)$ running through a basis of $H_{Y'(n/2)}$.  By Corollary
\ref{cor:BDbasis}, this implies that the image of $C$ is a basis of
$A$. Injectivity of $C$ follows from injectivity of $C_\theta$, $C_t$ $(0\le
t\le n/2)$, $C'_{n/2}$ if $n$ is even, and Theorem \ref{th:main}, which
guarantees that no collapses of dimensions of the individual parts occur.

\nl(C2).
Clearly, $*$ is an $S$-linear anti-involution.
Let $t\in\{0,\ldots,\lfloor n/2\rfloor\}$,
$\lambda\in\Lambda_t$, and $(B,x),(B',y)\in T(\lambda)$.
Then $(\rr_{B,n}C_t(x,y)\rr_{B',n}^\op)^\op
=\rr_{B',n}C_t(x,y)^\op \rr_{B,n}^\op$,  so, in order to establish
$\big(C((B,x),(B',y))\big)^*=C((B',y),(B,x))$,
it suffices to verify
that $C_t(x,y)^\op $ coincides with $C_t(y,x)$.
Now $*_t$ on
$H_Y{(t)}$ mod $J_{t+1}$ coincides with opposition, so
modulo $J_{t+1}$ we have $C_t(x,y)^\op =  C_t(x,y)^{*_t} = C_t(y,x)$
by the cellularity of $(\Lambda_t,T_t,C_t,*_t)$.
On the other hand, as $H_{Y(t)}$ is invariant under opposition,
see Corollary \ref{cor:HYt}(\ref{opp-invH}), and contains the values of $C_t$,
it contains $C_t(x,y)^\op - C_t(y,x)$, so
$C_t(x,y)^\op-C_t(y,x)\in H_{Y(t)}\cap J_{t+1} = \{0\}$, whence
$C_t(x,y)^\op=C_t(y,x)$, as required.

The case of $\lambda\in\Lambda'_{n/2}$ for $n$ even is similar.  If $\lambda
\in\Lambda_\theta$ and $x,y\in T(\lambda)$, then $C(x,y)^* = C(y,x)$ is
immediate from the cellularity of
$(\Lambda_\theta,T_\theta,C_\theta,*_\theta)$.

\nl(C3).  Let $\lambda\in\Lambda_t$ and $(B,x),(B',y)\in T(\lambda)$.
Fix $i\in \{1,\ldots,n\}$.
It clearly suffices to prove the formulas for $a$ running over the generators
$g_i$ and $e_i$ of $\BMW(\D_n)$.

By choice of $C_t$, we have
$C_t(x,y)\in H_{Y(t)}$, and,  see Proposition \ref{prop:hecke}(i),
$C_t(x,y)=\rho(\he_{Y(t)}) C_t(x,y)$.
According to Corollary~\ref{cor:HYt}(ii), there is $z_{B,i}\in H_{Y(t)}$,
depending only on $B$ and $i$,
such that $g_i\rr_{B,n}\rho(\he_{Y(t)}) \in \rr_{g_iB,n}z_{B,i} + J_{t+1}$.
As $(\Lambda_t,T_t,C_t,*_t)$ is a cell datum for
$H_{Y(t)}$ mod $J_{t+1}$,
for each $u\in T_t(\lambda)$,
there are $\nu_{i}(u,B,x)\in S$, independent of $B'$ and $y$,
such that
\begin{eqnarray*}
z_{B,i}C_t(x,y)&\in&\sum_{u \in T_t(\lambda)}
\nu_{i}(u,B,x)C_t(u,y) + (H_{Y(t)})_{<\lambda} + J_{t+1}.
\end{eqnarray*}
Since both  $(H_{Y(t)})_{<\lambda}$  and $J_{t+1}$ are contained in
$A_{<\lambda}$, we find
\begin{eqnarray*}
g_iC((B,x),(B',y)) & = &
g_i\rr_{B,n}\rho(\he_{Y(t)})C_t(x,y)\rr_{B',n}^\op\\
&\in&\rr_{r_iB,n}z_{B,i}C_t(x,y)\rr_{B',n}^\op + A_{<\lambda}\\
&=&\sum_{u\in T_t(\lambda)}
\nu_i(u,B,x)\rr_{r_iB,n}C_t(u,y)\rr_{B',n}^\op + A_{<\lambda}\\
&=&\sum_{u\in T_t(\lambda)}
\nu_i(u,B,x)C((r_iB,u),(B',y)) + A_{<\lambda},
\end{eqnarray*}
as required.

Rewriting (RSrr) to $e_i = lm^{-1}(g_i^2+mg_i-1)$, we see that, if $m^{-1}\in
S$, the proper behavior of the cell data under left multiplication by
$e_i$ is taken care of by the above formulae for $g_i$.  A proof in full
generality can be given that is
similar to the above proof for $g_i$ using Corollary~\ref{cor:HYt}(ii).

For
$\lambda\in\Lambda_\theta$, the formulas are straight from those for $\Theta'$
as $g_1a = g_2a$ and $e_1a= e_2a$ for each
$a\in\Theta'$.

This establishes that $(\Lambda,T,C,*)$ is a cell datum for $A$
and so completes the proof of Theorem~\ref{th:cellular}.

Alternatively, the information we have provided shows $\BMW(\D_n)$ is an
iterated inflation of Hecke algebras of type $\D_n$, $\D_{n-2t}$,
and $\A_{n-1-2s}$ for $s,t\ge1$ and so are cellular by \cite{KX}.

\section{Discussion}
\label{sec:discussion}

The following consequence of Corollary \ref{cor:finalreduction}
will be of use in \cite{CW4}.  It involves the
products $A_{B,n} = a_{B,n}\he_Y$ where $B\in WB_Y$.

\begin{Thm}\label{prop:aB}
Let $Y \in\MY$.
For each $B\in W B_Y$ there is, up to homogeneous
equivalence and powers of $\delta$, a unique word $A_{B,n}$ in $F_n\he_Y$
satisfying
the following three properties for each node $i$ of $\D_n$.
\begin{enumerate}[(i)]
\item
$r_i A_{B,n}\homog A_{r_iB,n}h$ for some $h\in U_Y$.
Furthermore, if $r_iB>B$,
then $h$ is the identity $\he_Y$ of $U_Y$.
\item If $|e_iB| = |B|$, then
$e_iA_{B,n}\homog A_{e_iB,n}h$ for some $h\in \delta^{\Z}U_Y$
and $\het(e_iB)\le \het(B)$.
\item If $|e_iB| > |B|$, then
$e_iA_{B,n}$ reduces to an element of $\BrM(\D_n) e_{U}\BrM(\D_n)$
for some set of nodes $U$ strictly containing $Y$.
\end{enumerate}
\label{th:aBcharacterization}
\end{Thm}

In \cite[Proposition~3.1]{CGW2} it is shown that there is a natural
order on each $W$-orbit in $\AO$, and in fact, \cite[Corollary~3.6]{CGW2},
each such orbit has a unique maximal element under this order.  The ordering
is also involved in a notion of height for elements of $\AO$, denoted
$\het(B)$ for $B\in\AO$, which satisfies $\het(B)<\het(C)$ whenever $B$ and
$C$ are in the same $W$-orbit in $\AO$ and satisfy $B<C$. Moreover, if
$r_iB>B$, then $\het(r_iB) = \het(B)+1$.  There are certain minimal elements
$Y\in \AO$ (such as the sets $B_Y$ for $Y\in\MY$ described above).  Then
$\het(B)$ will the distance  to the maximal element
in the Hasse diagram of the part of the poset
$WB=WB_Y$.  In particular,
$\het(B_Y)=0$.
The word $A_{B,n}$ has height $\het(B)$ and moves
$\emptyset$ to $B$ in the left action: $A_{B,n} (\emptyset) = B$.

The words $A_{B,n}$ are as given by Notation~\ref{not:ABN} using Lemma
\ref{lm:allthesame} and involve an ordering of the roots of $B$.  Height
considerations as above give an algorithm for choosing a representative for
$A_{B,n}$.  To begin, pick $\b_n$ to be a root of smallest height in $B$.
Then $\b_{n-2}$ should be a root of smallest height in $ a_{\b_n,n} e_n B
\setminus \{\a_n\}$ in case $B\in WB_{Y(t)}$.  Continue at each step picking
the next root $\b_{n-2s}$ as one of smallest height from the roots
remaining.  Similarly for other $Y\in\MY$.

An alternative proof of the results of this paper, using the methods of
\cite{CW4} is possible. The proof in that paper deals with the case $M=\E_n$
$(n=6,7,8)$ and involves a search of a finite number of finite posets.  The
search can be avoided in the case $M=\D_n$ by using the specific structure
of the root system and induction on $n$.  In \cite{CW4}, the definition of
$A_{B,n}$ (which is denoted $a_B$ there) is given by the following algorithm
which need not be the same as the one above.

\begin{Alg}\label{def:aB}
\rm
Given $B\in\AO$, determine a word $A_{B,n}$ of minimal height with
$A_{B,n}(\emptyset) = B$.

\begin{enumerate}[(i)]
\item
If the number of simple nodes in $B$ is $t$, then $A_{B,n}$ is a product of
$e_i$ which moves these simple nodes to $Y$ as described in \cite{CW4}.
\item
If $r_kB<B$ for some node $k$, then $A_{B,n} = r_kA_{r_kB,n}$.

\item
Otherwise, there are adjacent nodes $j,k$ with $\a_j\in B$; then $A_{B,n} =
e_jA_{e_kB,n}$ where $A_{e_kB,n}$ has been defined inductively.  Here $e_kB =
r_jr_kB$, $e_je_kB = B$, and $\het(e_kB) = \het(B)$.
\end{enumerate}
\end{Alg}

\np The main result of this paper concerns an upper bound for the BMW
algebra, given by a presentation. A lower bound, as can be seen in the proof
of the main theorem, is in \cite{CFW}. In particular, the current results
finish the proof of the main theorem in the paper \cite{CGW3} on the tangle
algebra $\KT(\D_n)$, which gives a topological depiction of $\BMW(\D_n)$.

On the level of the Brauer algebra, for a monomial $a$, the admissible set
$B = a(\emptyset)$ determines the connections of the horizontal strands at
the top in the following way: if $\eps_i-\eps_j$ belongs to $B$, then there
is a horizontal strand from top node $i$ to top node $j$ that does not go
around the pole.  If $\eps_i+\eps_j$ belongs to $B$, then there is a
horizontal strand from top node $i$ to top node $j$ that goes around the
pole.  If $\eps_i-\eps_j$ and $\eps_i+\eps_j$ both belong to $B$, then
$\Theta$, the pair of loops going around the pole as defined in \cite{CGW3},
belongs to the tangle.
In addition, the paper
\cite{CGW3}  gives an alternative proof of the lower bound on the rank
of $\BMW(\D_n)$.

Furthermore, the ideal $\Theta'$ in $\BMW(\D_n)$ has a nice interpretation in
the tangle algebra $\KT(\D_n)$, where they are ordinary tangles with no
loops around the pole and with coefficient $\Theta$, as described in
\cite{CGW3}. Here there must be at least one horizontal strand at the top and
one at the bottom. These tangles span
the ideal corresponding to $\mu^{-1}(Q)$, for $Q$ as
in Proposition \ref{lm:Theta1}. This ideal
is easily seen to be the ideal in $\BMW(\A_{n-1})$
identified in Section \ref{sec:conclusion}.

\end{document}